\documentclass[11pt]{amsart}
\usepackage[myheadings]{fullpage}
\usepackage{enumerate}
\usepackage{mathtools}
\usepackage{amssymb}
\usepackage{float}
\usepackage{xcolor}
\usepackage{hyperref}
\hypersetup{colorlinks=true, linkcolor=blue,citecolor=blue}

\theoremstyle{plain}
\newtheorem{theorem}{Theorem}[section]
\newtheorem{lemma}[theorem]{Lemma}
\newtheorem{prop}[theorem]{Proposition}
\newtheorem{corollary}[theorem]{Corollary}

\theoremstyle{definition}
\newtheorem{definition}[theorem]{Definition}
\newtheorem{remark}[theorem]{Remark}
\newtheorem{Example}[theorem]{Example}

\begin{document}
\title{Alexander-Markov correspondence for doodles on closed surfaces}
\author{Komal Negi}
\author{Mahender Singh}
	
\address{Department of Mathematical Sciences, Indian Institute of Science Education and Research (IISER) Mohali, Sector 81,  S. A. S. Nagar, P. O. Manauli, Punjab 140306, India.}
	\email{komalnegi@iisermohali.ac.in}
	\email{mahender@iisermohali.ac.in}

\date{}


\begin{abstract}
In this paper, we introduce twisted virtual doodles, defined as stable equivalence classes of immersed circles on closed surfaces that may be non-orientable. These objects admit planar representative diagrams, considered up to a suitable set of Reidemeister-type moves. To develop the associated braid-theoretic framework, we define twisted virtual twin groups as natural extensions of virtual twin groups, and establish Alexander- and Markov-type theorems in this set-up. This shows that twisted virtual doodles unify and extend both classical and virtual doodle theories. We further investigate the structure of the pure twisted virtual twin group, providing a presentation and deriving several structural and combinatorial properties. In particular, we obtain two interesting decompositions of the twisted virtual twin group and its pure subgroup, from which it follows that both groups have trivial center and are residually finite as well as Hopfian.
\end{abstract}

\maketitle


\section{Introduction}
Classically, isotopy classes of links in the 3-sphere (equivalently, in the thickened 2-sphere) correspond bijectively to Markov equivalence classes of Artin braids. Kauffman \cite{MR1721925} later generalized this framework by considering stable embeddings of disjoint unions of circles in thickened compact oriented surfaces via oriented Gauss codes, thereby introducing the theory of virtual links. When the surface being thickened is a sphere, one recovers classical links; however, more general surfaces give rise to new phenomena not present in the classical setting. Virtual braid groups were subsequently defined in \cite{MR2351010}, and analogues of the Alexander and Markov theorems were established, extending the classical braid-link correspondence to the virtual category. More recently, Bourgoin \cite{MR2443243} introduced twisted virtual links, defined as stable equivalence classes of oriented links in orientable 3-manifolds that are orientation  $I$-bundles over closed surfaces, including non-orientable ones. This theory simultaneously subsumes both virtual and projective links. Twisted virtual braid groups were later introduced by Negi, Prabhakar, and Kamada \cite{MR4812005}, who also established twisted virtual Alexander- and Markov-type theorems, thereby extending the braid-theoretic approach to this broader setting.
\par

As a two-dimensional counterpart of classical knot theory, the study of immersed circles on surfaces was initiated by Fenn and Taylor \cite{MR0547452}  in the case of the 2-sphere, and subsequently extended by Khovanov \cite{MR1370644} to arbitrary closed oriented surfaces. In this framework, the analogue of the braid group is given by a class of right-angled Coxeter groups known as twin groups (also referred to as planar braid groups). A virtual extension of this theory was developed by Bartholomew, Fenn, Kamada, and Kamada \cite{MR3876348}, leading to the introduction of the corresponding braid-like objects, termed virtual twin groups, by Bardakov, Singh, and Vesnin \cite{MR4027588}. Furthermore, Alexander- and Markov-type theorems for virtual doodles were established by Nanda and Singh \cite{MR4209535}, providing a foundational correspondence between virtual doodles and virtual twin groups. The combinatorial properties of these groups, including their automorphism groups, crystallographic quotients, simplicial structures, and twisted conjugacy classes, have since been studied in detail in a series of works \cite{MR4805520, MR4607569, MR4750162, MR4685104, MR4651964, MR4739232}.
\par

In this paper, we introduce the notion of twisted virtual doodles, defined as stable equivalence classes of doodles on closed surfaces that may be non-orientable. Analogous to other knot and link theories, such objects admit planar representatives, called twisted virtual doodle diagrams, considered up to an appropriate set of Reidemeister-type moves. To develop the corresponding braid-theoretic framework, we define twisted virtual twin groups as extensions of virtual twin groups, and establish Alexander- and Markov-type theorems in this setting. In this way, twisted virtual doodles naturally extend and unify both classical and virtual doodles. Furthermore, we examine the structure of the pure twisted virtual twin group and derive several of its algebraic and combinatorial properties. These correspondences can be  summarised as follows:
\begin{eqnarray*}
\textstyle\bigcup_{n \ge 2} T_n /_{\text{Markov equivalence}} &\longleftrightarrow & \text{Equivalence classes of doodles on the 2-sphere},\\
\textstyle\bigcup_{n \ge 2} VT_n /_{\text{Markov equivalence}}  &\longleftrightarrow & \text{Stable equivalence classes of doodles on closed}\\
& & \text{oriented surfaces},\\
\textstyle\bigcup_{n \ge 2} TVT_n /_{\text{Markov equivalence}}&\longleftrightarrow & \text{Stable equivalence classes of doodles on closed }\\
&  & \text{surfaces not necessarily oriented}.
\end{eqnarray*}

This paper is organized as follows. In Section \ref{ds}, we introduce twisted virtual doodles as stable equivalence classes of immersed circles on closed surfaces, allowing the surface to be non-orientable, and define twisted virtual doodle diagrams on the plane. The main result of this section (Theorem \ref{thm bijection}) establishes a bijection between stable equivalence classes of twisted virtual doodles on surfaces and equivalence classes of twisted virtual doodle diagrams on the plane. In Section ~\ref{sect:braid}, to develop the corresponding braid-theoretic framework, we introduce twisted virtual twins as equivalence classes of certain configurations of polygonal curves in the strip $\mathbb{R} \times [0, 1]$ and give an explicit presentation of the resulting group $TVT_n$, called the twisted virtual twin group on $n$ strands. This group turns out to be a natural extension of the virtual twin group $VT_n$. In Section \ref{sect:Alexander}, we prove an Alexander-type theorem (Theorem \ref{theoremAlexanderB}), asserting that every twisted virtual doodle is equivalent to the closure of a twisted virtual twin diagram. In Section \ref{sect:Markov}, we introduce Markov-type moves on twisted virtual twins and prove that two twisted virtual twin diagrams have equivalent closures (as twisted virtual doodles) if and only if they are Markov equivalent (Theorem \ref{theorem:MarkovA}). As in other braid theories, there is a natural surjection $TVT_n \to S_n$, whose kernel is the pure twisted virtual twin group $PTVT_n$. In Section \ref{sec pure}, we provide a presentation of this subgroup (Theorem \ref{sTVP_n}). Finally, in Section \ref{sec decompositions}, we obtain two decompositions of $PTVT_n$.  The first (Theorems \ref{decom1}) expresses it as an iterated semidirect product of free products of infinite-rank free groups and an involution. The second (Theorem \ref{decom2}) describes it as a semidirect product of an irreducible right-angled Artin group with an elementary abelian 2-group. As a consequence, we show that both $TVT_n$ and $PTVT_n$ have trivial center and are residually finite as well as Hopfian for all $n\ge 2$ (Corollaries \ref{cor trivial center} and  \ref{cor:3.5}).

\medskip

\section{Twisted virtual doodles and their planar diagrams}\label{ds}

\subsection{Twisted virtual doodles on surfaces} We begin by introducing the main objects of our study.

\begin{definition}
A {\it doodle} $D$ on a closed surface  $\Sigma$ (which may be orientable or non-orientable) is the image of a smooth immersion 
$$f : \bigsqcup_{i=1}^n \mathbb{S}^1 \longrightarrow \Sigma$$
of disjoint union of circles such that
$$|f^{-1}(f(x))| < 3 \quad \text{for all } x \in\bigsqcup_{i=1}^n \mathbb{S}^1.$$
That is, $D$ does not admit any triple or higher-order intersection points. The integer $n$ is called the {\it number of components} of $D$.
\end{definition}

We say that two doodles are {\it stably equivalent} if one can be transformed into the other by a finite sequence of following equivalences:

\begin{itemize}
    \item \textbf{Homeomorphic equivalence.}
Two doodles $D$ and $D'$ on surfaces $\Sigma$ and $\Sigma'$, respectively, are said to be homeomorphically equivalent if there exist a homeomorphism $ S : \Sigma \;\longrightarrow\; \Sigma'$ such that $S(D)=D'$.

\item \textbf{Homotopic equivalence.} Two doodles $D$ and $D'$ on a surface $\Sigma$ are said to be homotopically equivalent if one can be transformed into the other by a finite sequence of local moves on $\Sigma$ as illustrated in Figure~\ref{he}. These moves correspond to the creation or elimination of monogons and bigons.

\begin{figure}[ht]
  \centering
    \includegraphics[width=6cm]{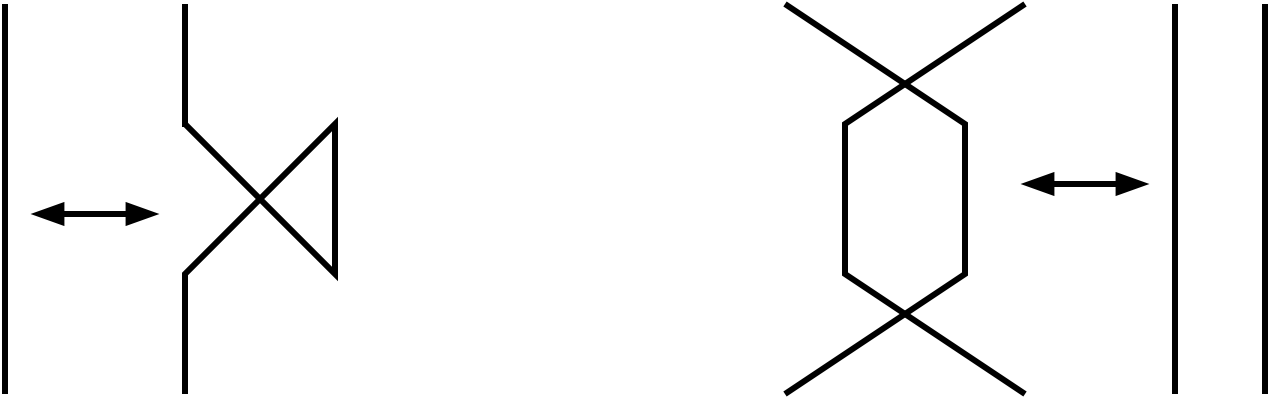}
        \caption{Moves for homotopic equivalence.}
        \label{he}
        \end{figure}
        
\item \textbf{Surgery equivalence.}
Surface surgery consists of a finite sequence of additions and removals of handles or cross caps that are disjoint from the doodle diagram. To perform a handle addition,   consider two closed discs on the surface that are disjoint from the doodle.  Remove the interiors of these discs and replace them with an annulus 
  \(\mathbb{S}^1 \times [0,1]\), attaching its boundary circles to the boundaries of the 
  two removed discs (see Figure~\ref{ss} A). 
  To perform a cross-cap addition, consider one closed disc on the surface that is disjoint from the doodle.  Remove the interior of the disc and replace it with a cross-cap, attaching its boundary circle to the boundary of the removed disc (see Figure~\ref{ss} B).
  
  In the reverse direction, we cut the surface along a non-separating, two-sided simple closed curve which is disjoint from the doodle and then cap off the two resultant boundary components with 2-discs. Here, a 2-sided simple closed curve is one that has a regular neighborhood homeomorphic to an annulus.

  \begin{figure}[ht]
  \centering
    \includegraphics[width=8cm]{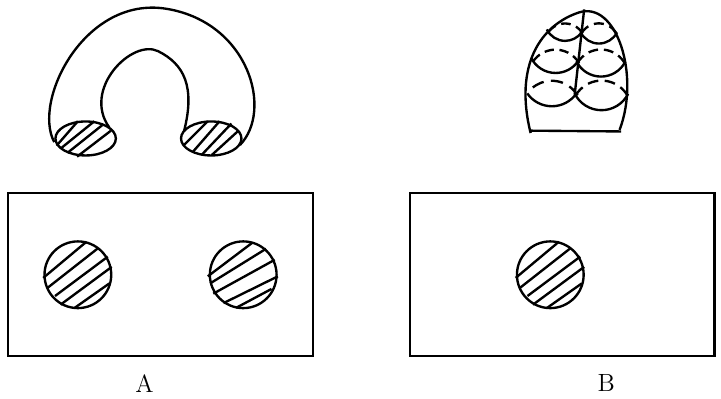}
        \caption{Surface surgery.}
        \label{ss}
        \end{figure}
\par

Furthermore, if the surface $\Sigma$ has a connected component $\Sigma_0$ that is disjoint   from the doodle $D$, then we simply ignore $\Sigma_0$ and consider the doodle $D$ on the surface $\Sigma \setminus \Sigma_0$.
\end{itemize}

\begin{definition}
The stable equivalence class of a doodle is called a  {\it  twisted virtual doodle}. Further, a twisted  virtual  doodle is called {\it trivial} if its stable equivalence class contains a representative that is a multicurve on the 2-sphere. Here, a multicurve is a finite collection of pairwise disjoint simple closed curves.

Henceforth, by abuse of terminology, we shall refer to a representative doodle as a twisted virtual doodle.

\end{definition}
\begin{Example}
Some trivial twisted virtual doodles on different surfaces are illustrated in Figure~\ref{trid}. In (A), the twisted virtual doodle is drawn on the torus $T$. By performing a surface surgery along a non-separating, two-sided curve, we see that it is equivalent to a twisted virtual doodle on the 2-sphere, which is a simple closed curve.  In (B), it is drawn on $T \# \mathbb{RP}^2$, the connected sum of a torus and a projective plane. By performing a surface surgery along a non-separating, two-sided curve, we see that it is equivalent to previous twisted virtual doodle. In (C), it is drawn on the Klein bottle $\mathbb{RP}^2 \# \mathbb{RP}^2$. By performing a surface surgery along a non-separating, two-sided curve which passes through both the crosscaps, we can show that it is equivalent to a twisted virtual doodle on the 2-sphere.

\begin{figure}[ht]
  \centering
    \includegraphics[width=14cm]{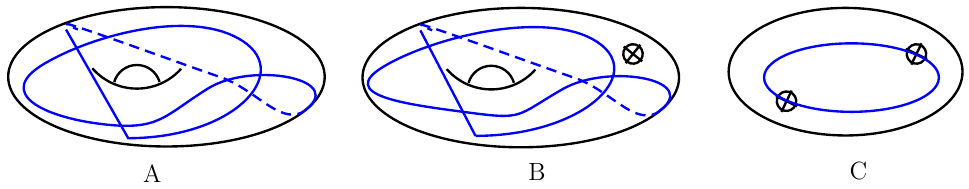}
        \caption{Examples of trivial twisted virtual doodles on surfaces.}
        \label{trid}
        \end{figure}
\end{Example}

\begin{Example}\label{non-trivial example}
Figure~\ref{ntrid} shows a non-trivial twisted virtual doodle on the surface $T \# T \# \mathbb{RP}^2$. A justification for its non-triviality is given in Remark \ref{remark justification}.

\begin{figure}[ht]
  \centering
    \includegraphics[width=5cm]{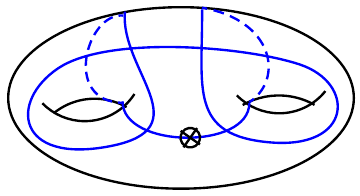}
        \caption{Example of a non-trivial twisted virtual doodle on a surface.}
        \label{ntrid}
        \end{figure}
\end{Example}
\medskip

\subsection{Twisted virtual doodle diagrams on the plane} In this subsection, we introduce diagrammatic representations of twisted virtual doodles.

\begin{definition}
A {\it twisted virtual doodle diagram} on the plane is a generic immersion of finitely many circles, admitting two type of crossings, namely,   the real crossings and the virtual crossings which are decorated with small circles. Furthermore, some of the edges are decorated with one or more bar-like symbols. Such diagrams can be regarded as 4-valent graphs on the plane with two types of vertices (classical and virtual), and two types of edges (with and without bars). A twisted  virtual doodle diagram is said to be {\it oriented} if an orientation is fixed on each of its components.
\end{definition}

\begin{figure}[ht]
  \centering
    \includegraphics[width=10cm]{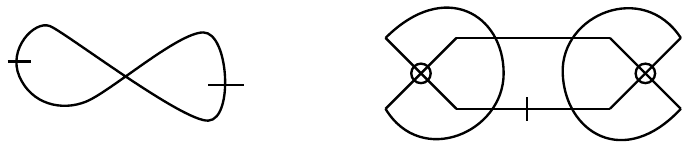}
        \caption{Examples of twisted virtual doodle diagrams.}    
        \label{extd}
        \end{figure}

See Figure~\ref{extd} for some examples.
 
\begin{figure}[ht]
  \centering
    \includegraphics[width=10cm]{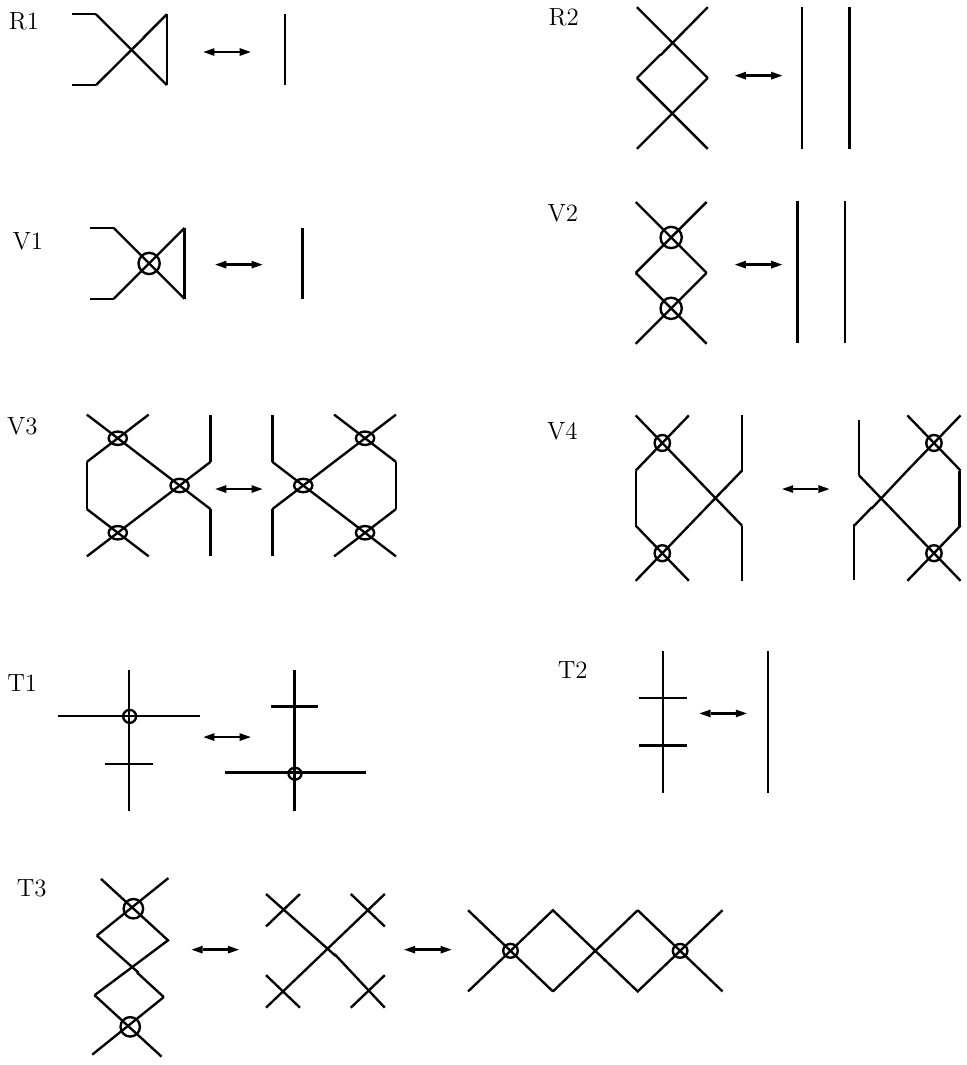}
        \caption{Extended Reidemeister moves for twisted virtual doodle diagrams.}
        \label{tm}
        \end{figure}

\begin{definition}
Two twisted virtual doodle diagrams on the plane are said to be {\it equivalent} if one can be transformed into the other by a finite sequence of R1, R2, V1, V2, V3, V4, T1, T2, and T3 moves, together with planar isotopies, where the moves are as shown in Figure~\ref{tm}.
\end{definition}

We will refer to an equivalence class of twisted virtual doodle diagrams simply as a twisted virtual doodle diagram.

\begin{definition}
Let $D$ be a twisted virtual doodle diagram. The \emph{parity of bars} of $D$ is defined as the parity of the total number of bars appearing in $D$, that is, whether the number of bars is even or odd.
\end{definition}

\begin{prop}
The parity of bars is an invariant of twisted virtual doodle diagrams. In other words, if two twisted virtual doodle diagrams are equivalent, then they admit the same parity of bars.
\end{prop}

\begin{proof}
The moves defining the equivalence of twisted virtual doodle diagrams either preserve the number of bars or change it by an even number. Hence, the parity of the number of bars remains unchanged under these moves, and therefore it is an invariant of twisted virtual doodle diagrams.
\end{proof}

As an example, Figure~\ref{extd} represents inequivalent twisted virtual doodle diagrams.

\begin{definition}
A twisted virtual doodle diagram is said to be  {\it trivial} if it is equivalent to a disjoint union of circles on the plane.
\end{definition}

\begin{remark}
The parity of bars yields a natural classification of twisted virtual doodle diagrams into the following two distinct classes:
\begin{enumerate}
    \item Diagrams with an even number of bars, which, in particular, includes the class of trivial doodles as well as virtual doodles (since they contain no bars).
    \item Diagrams with an odd number of bars.
\end{enumerate}
Thus, the parity of bars establishes the existence of non-trivial twisted virtual doodle diagrams which are not virtual doodle diagrams.
\end{remark}
\medskip

\begin{figure}[ht]
  \centering
    \includegraphics[width=12cm]{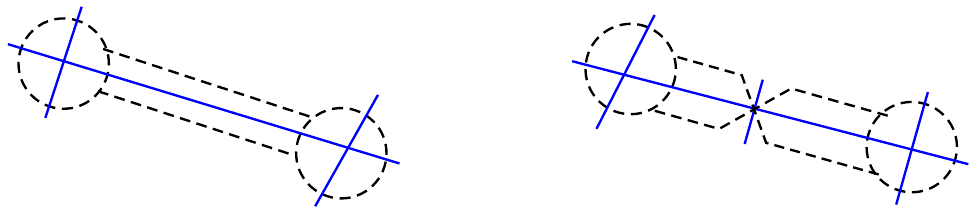}
        \caption{Attaching a band and a twisted band.}     
        \label{band}
        \end{figure}

We shall use the two operations shown in Figure~\ref{band} in the subsequent theorem. The operations represent attaching a band and a twisted band to a pair of two disks on the plane.

\begin{theorem}\label{thm bijection}
There is a bijection between the set of stable equivalence classes of twisted virtual doodles on surfaces and the set of equivalence classes of twisted virtual doodle diagrams.    
\end{theorem}

\begin{proof}
Let $D$ be a twisted virtual doodle on a closed surface \(\Sigma\). Let \(N(D)\) denote a  regular neighborhood of \(D\) in \(\Sigma\) such that $D$ is a deformation retract of $N(D)$. Consider the regular projection of \(N(D)\) onto the plane (we can think of it as flattening the \(N(D)\) onto the plane) and let $\widetilde{D}$ be the image of $D$ under the regular projection. We mark the projection of each crossing of $D$ as a real crossing on $\widetilde{D}$. Further, if two bands lie one above the other in $\Sigma$, then their projection onto the plane creates an intersection point on $\widetilde{D}$, which we mark as a virtual crossing. Finally, if a band on $\Sigma$ is a twisted band, then we indicate it by placing a bar on the corresponding arc on $\widetilde{D}$. Hence, we obtain a twisted virtual doodle diagram $\widetilde{D}$ from a  twisted virtual doodle $D$ on the surface $\Sigma$. It is straightforward to check that if two twisted virtual doodles are equivalent, then their corresponding twisted virtual doodle diagrams are also equivalent. Thus we have a well-defined map which sends the class of a twisted virtual doodle to the class of the corresponding twisted virtual doodle diagram on the plane.
\par
Let $\widetilde{D}$ be a twisted virtual doodle diagram on the plane. First, we construct a regular neighbourbood $N(\widetilde{D})$ of $\widetilde{D}$, which is a compact surface. We place a small disk around each real crossing of $\widetilde{D}$, and then connect these disks with bands or twisted bands along the arcs of $\widetilde{D}$. When an arc of  $\widetilde{D}$ does not contain a bar, then we join the corresponding disks with a  band.  When an arc of  $\widetilde{D}$ contains a bar, then we join the corresponding disks with a twisted  band. For each virtual crossing of  $\widetilde{D}$, we attach two bands or twisted bands stacked one above the other, depending on which arc has a bar. In this manner, the diagram  $\widetilde{D}$ thickens into a compact surface $N(\widetilde{D})$ having  $\widetilde{D}$ as its deformation retract.  Let $E$ be a compact surface such that
\(\partial N(\widetilde{D})\) and \(\partial E\) have the same number of connected components. Take a homeomorphism $f : \partial E \to \partial N(\widetilde{D})$ and consider the closed surface $\Sigma$ obtained by gluing $E$ and $N(\widetilde{D})$ along  their boundaries via $f$. The image of $\widetilde{D}$ in the surface $\Sigma$ gives a twisted virtual doodle $D$ on $\Sigma$. Let  $E'$ be another surface and $\Sigma'$ be obtained by gluing $E'$ and $\partial N(\widetilde{D})$ as above.  If $D'$ is the image of $\widetilde{D}$ in $\Sigma'$, then  applying surface surgeries show that $(D, \Sigma)$ and  $(D', \Sigma')$ are equivalent as twisted virtual doodles. Similarly, if $\widetilde{D}_1$ and $\widetilde{D}_2$ are equivalent twisted virtual doodle diagrams, then their corresponding twisted virtual doodles $(D_1, \Sigma_1)$ and $(D_2, \Sigma_2)$ are equivalent. This gives a well-defined map sending the class of a twisted virtual doodle diagram to the class of a twisted virtual doodle on a surface, and it is the inverse of the preceding map.
 \end{proof}

 \begin{remark}\label{remark justification}
 Theorem~\ref{thm bijection} allows us to conclude that Example \ref{non-trivial example} is a non-trivial twisted virtual doodle, as detected by the parity of bars.
 \end{remark}
\medskip

\section{Twisted virtual twin groups}\label{sect:braid}

\begin{definition}
A {\it twisted virtual twin diagram} on $n \ge 1$ strands is a union of $n$ smooth or polygonal curves in $\mathbb{R} \times [0,1]$, called {\it strands}, connecting the points $\{ (1,1), \ldots, (n,1) \}$ with the points $\{ (1,0), \ldots, (n, 0) \}$ in a one-to-one fasion such that the following conditions hold:
\begin{enumerate}
\item The stands are monotonic with respect to the second coordinate.
\item The stands intersect at  transverse double points, which are equipped with the labelling of being a real or a virtual crossing.
\item The stands may be labelled with {\it bars} which are disjoint from crossings, which are short arcs intersecting the strands transversely.
\end{enumerate}
\end{definition}
\begin{figure}[ht]
  \centering
    \includegraphics[width=3cm]{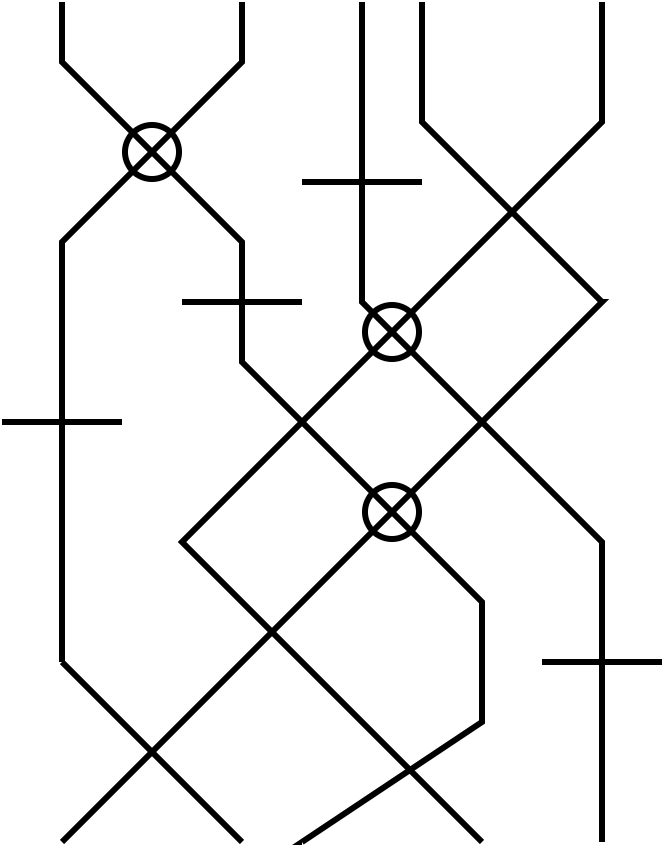}
        \caption{A twisted virtual twin diagram on five strands.}
        \label{exa}
        \end{figure}

See Figure~\ref{exa} for a twisted virtual twin diagram on five strands. Let $X(b)$ denote the set of crossings of a  twisted virtual twin diagram $b$ and the points on the strands where the bars intersect with the strands.  We say that the  twisted virtual twin diagram $b$ is  {\it faithful} if the projection onto the second coordinate $X(b) \to [0, 1]$ is injective. Note that, each twisted virtual twin diagram can be turned faithful by using planar isotopies.
\par

Next, we define three types of local moves on twisted virtual twin diagrams, which we refer to as {\it twisted Reidemeister twin moves}.  

\begin{itemize}
\item Figure~\ref{bmoves} shows the classical twin move R2, which is the planar analogue of the second Reidemeister move. 
\begin{figure}[ht]
  \centering
    \includegraphics[width=7.5cm]{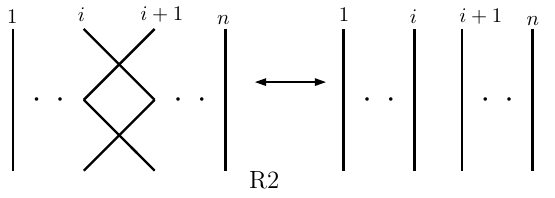}
        \caption{Classical twin move.}
        \label{bmoves}
        \end{figure}  
\item Figure~\ref{vbmoves} shows the virtual twin moves V2, V3 and V4,  which are the planar analogues of the virtual Reidemeister moves.
\begin{figure}[ht]
  \centering
    \includegraphics[width=15cm]{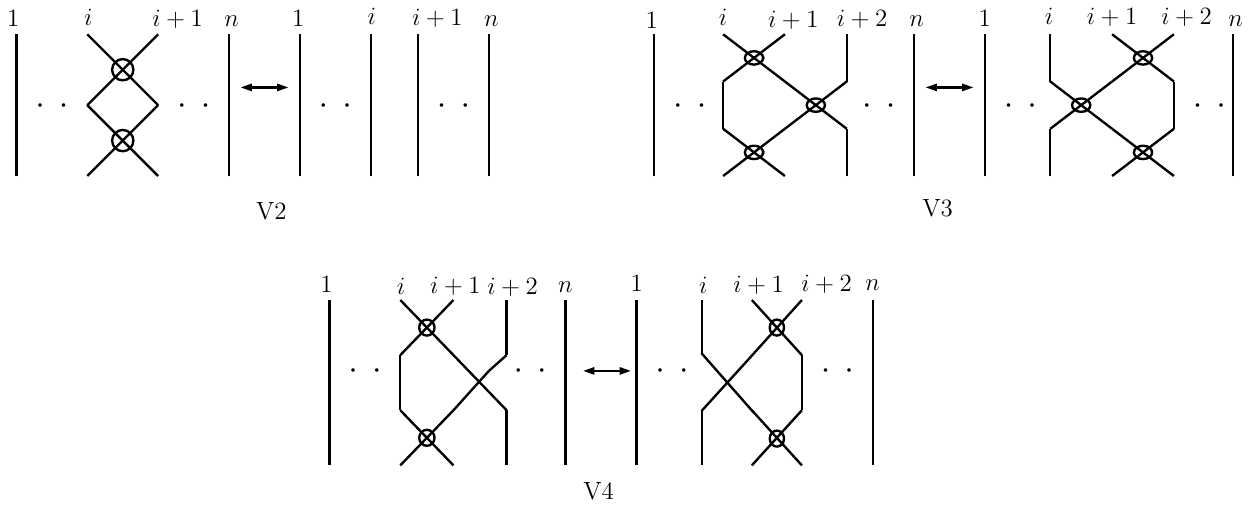}
        \caption{Virtual twin moves.}
        \label{vbmoves}
        \end{figure}  
\item Figure \ref{moves} shows the twisted twin moves T1, T2 and T3, which are the planar analogues of twisted Reidemeister moves.
 \begin{figure}[ht]
  \centering
    \includegraphics[width=14cm]{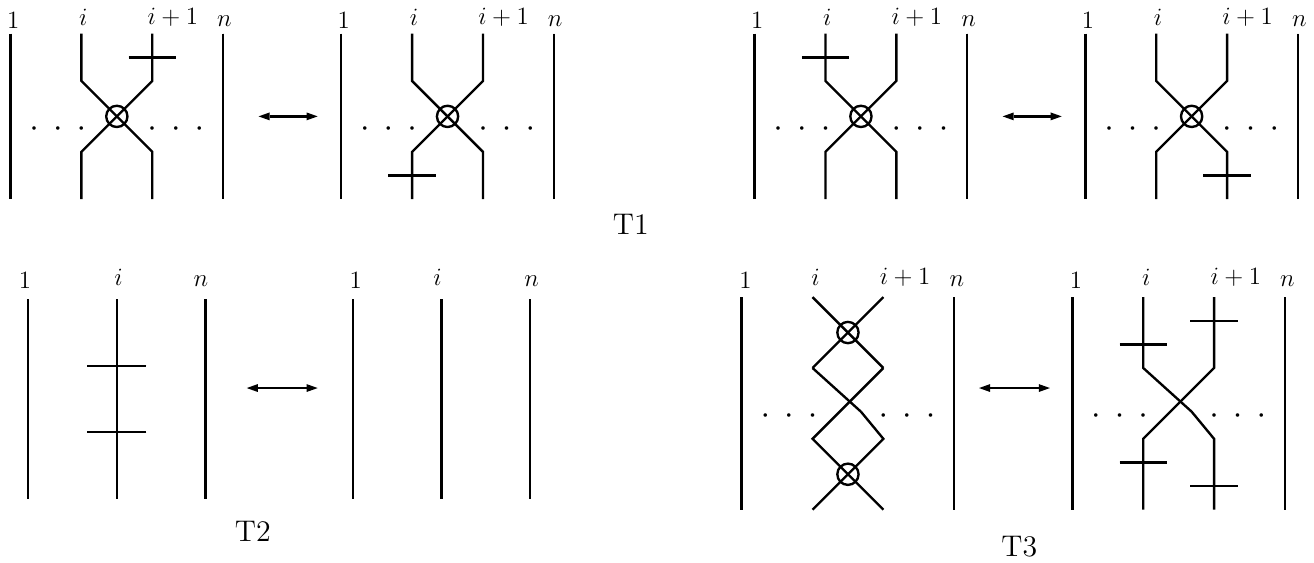}
        \caption{Twisted twin moves.}
        \label{moves}
        \end{figure}  
\end{itemize}

The preceding moves lead to the following definition.

\begin{definition}
Two twisted virtual twin diagrams $b$ and $b'$ on $n$ strands are said to be {\it equivalent} if there is a finite sequence $b_0, b_1, \dots, b_m$ of twisted virtual twin diagrams on $n$ strands
 with $b=b_0$ and $b'=b_m$ such that each $b_j$ is obtained from $b_{j-1}$ by one of the following: 
\begin{itemize}
    \item A planar isotopy satisfying the conditions (1)--(3) of a twisted virtual twin diagram.  
    \item A twisted Reidemeister twin move.  
\end{itemize}   
\end{definition}

\begin{definition}
      A {\it twisted virtual twin on $n$ strands} is an equivalence class of twisted virtual twin diagrams on $n$ strands.  
\end{definition}

For the sake of convenience, we denote a twisted virtual twin and a twisted virtual twin diagram representing it by  the same notation. The set of all twisted virtual twins on $n$ strands forms a group $\widetilde{TVT}_n$, where the product $b \, b'$ is defined by placing $b$ over $b'$ and then reparametrizing the interval back to $[0,1]$, similar to  products defined in other braid theories. 
\par

\begin{figure}[ht]
  \centering
    \includegraphics[width=12cm]{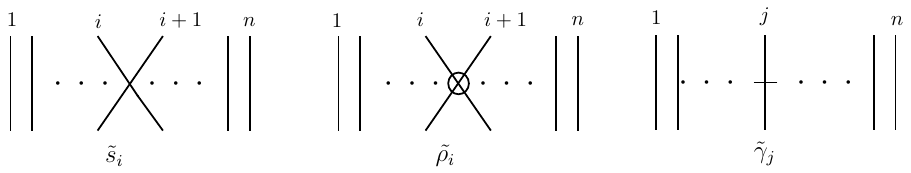}
        \caption{Generators for the twisted virtual twin group.}
        \label{gen}
        \end{figure}

Next, we give a generating set for $\widetilde{TVT}_n$. Let $\tilde{s}_i$, $\tilde{\rho}_i$ $(1 \le i \le n-1)$ and $\tilde{\gamma}_j$ $(1 \le j \le n)$ be 
twisted virtual twin diagrams as shown in Figure~\ref{gen}. By representing each twisted virtual twin by a faithful twisted virtual twin diagram, we see that $\widetilde{TVT}_n$ is generated by the set 
$$\{\tilde{s}_i, \tilde{\rho}_i, \tilde{\gamma}_j \mid 1 \le i \le n-1 ~~ \textrm{and} ~~ 1 \le j \le  n \}.$$

For each $n \ge 1$, let $TVT_n$ be an abstract group generated by the set  
$$\{s_i, \rho_i, \gamma_j \mid 1 \le i \le n-1 ~~ \textrm{and} ~~ 1 \le j \le  n \}$$
subject to the following defining relations:
    \begin{align}
        s_i^2 & = 1  & \text{ for } & 1\le i \le n-1,\label{rel-height-s2}\\
        s_i s_j & = s_j s_i  & \text{ for } & |i-j| > 1, \label{rel-height-ss}\\
        \rho_i^2 & = 1  & \text{ for } & 1\le i \le n-1, \label{rel-inverse-v}\\
        \rho_i \rho_j & = \rho_j \rho_i & \text{ for } & |i-j| > 1, \label{rel-height-vv}\\
        \rho_i \rho_{i+1} \rho_i & = \rho_{i+1} \rho_i \rho_{i+1} & \text{ for } & 1 \le i \le n-2, \label{rel-vvv}\\
        s_i \rho_j & = \rho_j s_i &  \text{ for } & |i-j| >1, \label{rel-height-sv}\\
        \rho_i s_{i+1} \rho_i & = \rho_{i+1} s_i \rho_{i+1} & \text{ for } & 1 \le i \le  n-2, \label{rel-vsv}\\
        \gamma_i^2 & = 1 & \text{ for } & 1\le i \le n, \label{rel-inverse-b}\\  
        \gamma_i \gamma_j & = \gamma_j \gamma_i   & \text{ for } & 1\le i,j \le n, \label{rel-height-bb} \\
        \gamma_j \rho_i & = \rho_i \gamma_j & \text{ for } & j\neq i, i+1, \label{rel-height-bv}\\
        s_i\gamma_j & = \gamma_js_i & \text{ for } & j\neq i, i+1, \label{rel-height-sb}\\
        \gamma_{i+1} \rho_i & = \rho_{i} \gamma_i & \text{ for } & 1\le i \le n-1, \label{rel-bv} \\
        \rho_{i} s_i \rho_{i} & = \gamma_{i+1} \gamma_i s_{i} \gamma_i \gamma_{i+1} & \text{ for } &  1\le i \le n-1. \label{rel-twist-III}
    \end{align}

\begin{theorem}\label{thm:StandardPresentation}
For each $n \ge 1$, we have $TVT_n \cong \widetilde{TVT}_n$. 
\end{theorem}

\begin{proof}
It follows from the definition of equivalence of two twisted virtual twin diagrams on $ n$ strands that the generators $\tilde{s_i}$, $\tilde{\rho_i}$ and $\tilde{\gamma_i}$ satisfy the following relations:
   \begin{align*}
        \tilde{s_i}^2 & = 1  & \text{ for } & 1\le i \le n-1,\\
        \tilde{s_i} \tilde{ s_j} & = \tilde{s_j} \tilde{ s_i}  & \text{ for } & |i-j| > 1,\\
        \tilde{\rho_i}^2 & = 1  & \text{ for } & 1\le i \le n-1,\\
        \tilde{\rho_i} \tilde{ \rho_j} & = \tilde{\rho_j} \tilde{\rho_i} & \text{ for } & |i-j| > 1,\\
        \tilde{\rho_i} \tilde{\rho_{i+1}} \tilde{\rho_i} & = \tilde{\rho_{i+1}} \tilde{\rho_i} \tilde{\rho_{i+1}} & \text{ for } & 1\le i \le n-2,\\
        \tilde{s_i} \tilde{ \rho_j} & = \tilde{\rho_j} \tilde{s_i} &  \text{ for } & |i-j| >1,\\
        \tilde{\rho_i} \tilde{s_{i+1}} \tilde{\rho_i} & = \tilde{\rho_{i+1}} \tilde{s_i} \tilde{\rho_{i+1}} & \text{ for } & 1\le i \le n-2,\\
        \tilde{\gamma_i}^2 & = 1 & \text{ for } & 1\le i \le n,\\  
        \tilde{\gamma_i} \tilde{\gamma_j} & = \tilde{\gamma_j} \tilde{\gamma_i}   & \text{ for } & 1\le i,j \le n, \\
        \tilde{\gamma_j} \tilde{\rho_i} & = \tilde{\rho_i} \tilde{\gamma_j} & \text{ for } & j\neq i, i+1,\\
        \tilde{s_i}\tilde{\gamma_j} & = \tilde{\gamma_j}\tilde{s_i} & \text{ for } & j\neq i, i+1,\\
        \tilde{\gamma_{i+1}} \tilde{\rho_i} & = \tilde{\rho_{i}} \tilde{\gamma_i} & \text{ for } & 1\le i \le n-1, \\
       \tilde{ \rho_{i}} \tilde{s_i} \tilde{\rho_{i}} & = \tilde{\gamma_{i+1}} \tilde{\gamma_i} \tilde{s_{i}} \tilde{\gamma_i} \tilde{\gamma_{i+1}} & \text{ for } &  1\le i \le n-1. 
    \end{align*}
Hence, the map $\psi_n: TVT_n \to \widetilde{TVT_n}$ defined on the generators by
$$\psi_n(s_i)=\tilde{s_i}, \quad \psi_n(\rho_i)=\tilde{\rho_i} \quad \textrm{and} \quad \psi_n(\gamma_j)=\tilde{\gamma_j}$$
for $1 \le i \le n-1$ and $1 \le j \le n$, extends to a group homomorphism. Clearly, the map $\psi_n$ is surjective since $\widetilde{TVT_n}$ is generated by $\{\tilde{s}_i, \tilde{\rho}_i, \tilde{\gamma}_j \mid 1 \le i \le n-1 ~~ \textrm{and} ~~ 1 \le j \le  n \}$. Next, we define $$\phi_n: \widetilde{TVT_n} \to TVT_n$$ 
by $\phi_n(\tilde{w})=w$, where $w$ is the word obtained from $\tilde{w}$ by replacing $\tilde{s_i}, \tilde{\rho_i}, \tilde{\gamma_j}$ by $s_i, \rho_i, \gamma_j$, respectively. Using the equivalence of twisted virtual twin diagrams, we can verify that $\phi_n$ is indeed well-defined. Since $\phi_n \circ \psi_n$ is the identity map, it follows that $\psi_n$ is injective, and hence an isomorphism.
\end{proof}

\begin{definition}
For each $n \ge 1$, the group $TVT_n$ (or $\widetilde{TVT}_n$) is called the {\it twisted virtual twin group on $n$ strands}.
\end{definition}

\begin{remark} 
There are two kinds of twisted virtual twin moves of type T1 as shown in Figure~\ref{moves}. The left one corresponds to the relation $\eqref{rel-bv}$, while the right one corresponds to the following relation
    \begin{equation} \label{rel-vb}
       \gamma_i \rho_i   = \rho_i \gamma_{i+1} ~~ \text{ for } ~~ 1 \le i \le n-1. 
     \end{equation}
Using $\eqref{rel-inverse-v}$, we see that the relation $\eqref{rel-bv}$ is equivalent to $\eqref{rel-vb}$. Thus, we may replace $\eqref{rel-bv}$ by $\eqref{rel-vb}$ in the presentation of $TVT_n$.
\end{remark}
\medskip


\section{Alexander theorem for twisted virtual doodles}\label{sect:Alexander}

Given an oriented twisted virtual doodle diagram $D$, we set the following conventions:
\begin{itemize}
     \item Let $\mu(D)$ be the number of components of $D$.
    \item Let $V_R(D)$ be the set of all real crossings of $D$.
    \item Let $B(D)$ be the set of all bars in $D$.
    \item Let $N(c)$ be a regular neighbourhood of $c$ for each $c \in V_R(D) \cup B(D)$.
     \begin{figure}[ht]
  \centering
    \includegraphics[width=3cm]{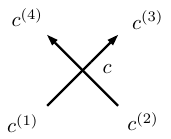}
        \caption{Boundary points of $N(c) \cap D$.}
        \label{r}
        \end{figure}
    \item For each $c \in V_R(D)$, let $c^{(1)}, c^{(2)}, c^{(3)}$, and $c^{(4)}$ denote the four points of $\partial N(c) \cap D$  (see Figure~\ref{r}). 
   
     \item For each $\gamma \in B(D)$, let $\gamma^{(1)}$ and $\gamma^{(2)}$ denote the two points of $\partial N(\gamma) \cap D$ (see Figure \ref{b}).
      \begin{figure}[ht]
  \centering
    \includegraphics[width=1cm]{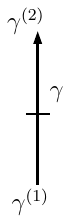}
        \caption{Boundary points of $N(\gamma) \cap D$.}\label{b}
        \end{figure}
     \item Let $W(D)=\overline{\mathbb{R}^2 \setminus \cup_{v\in V_R(D) \cup B(D)} N(v)}$, the closure of the set.
     \item Let $V^{\partial}_R(D) = \{c^{(j)} \mid c \in V_R(D),~ 1\leq j \leq 4\}$ and $B^{\partial}(D) = \{\gamma^{(j)} \mid \gamma \in B(D), ~1\leq j \leq 2\}.$ 
     \item Let $D  \cap W(D)$ be the intersection of $D$ and $W(D)$, which is a union of some oriented arcs generically immersed in $W(D)$ such that the double points are virtual crossings of $D$, and the set of boundary points of these arcs is $V^{\partial}_R(D) \cup B^{\partial}(D).$
     \item Let $X(D)$ be the subset of $(V^{\partial}_R(D) \cup B^{\partial}(D))^2$ consisting of tuples $(a,b)$ such that $D \cap W(D)$ admits an arc that starts from $a$ and ends at $b$. Further, the arc is allowed to have a virtual crossing.
\end{itemize}

\begin{definition}
    The {\it Gauss data} of an oriented twisted virtual doodle diagram $D$ is the quadruple $$\big(V_R(D), B(D), X(D), \mu(D)\big).$$
\end{definition}

\begin{Example}
Let $D$ be the oriented twisted virtual doodle diagram as in Figure~\ref{gd}. Then the real crossings of $D$ are $c_1$ and $c_2$, and the Gauss data is  
$$\Bigg(\{c_1,c_2\}, \{\gamma_1\}, \big\{ (c_1^{(4)}, c_2^{(2)}), (c_2^{(3)}, c_1^{(2)}), (c_2^{(4)}, \gamma_1^{(1)}), (\gamma_1^{(2)}, c_1^{(1)}), (c_1^{(3)}, c_2^{(1)}) \big\},1 \Bigg).$$
 
\begin{figure}[ht]
  \centering
    \includegraphics[width=4cm]{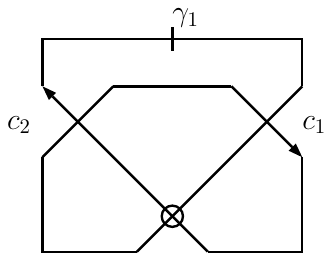}
        \caption{An oriented twisted virtual doodle diagram.}
        \label{gd}
        \end{figure}
\end{Example}

\begin{definition}
Two oriented twisted virtual doodle diagrams $D$ and $D'$ are said to have the {\it same Gauss data} if $\mu (D) = \mu (D')$ and there exists a bijection $g:V_R(D) \cup B(D) \to V_R(D') \cup B(D')$ satisfying the following conditions:
\begin{enumerate}
    \item $g(V_R(D)) = V_R(D')$ and $g(B(D))= B(D')$.    
    \item $(a,b) \in X(D)$ if and only if $(g^{\partial}(a), g^{\partial}(b)) \in X(D')$, where  $g^{\partial}: V^{\partial}_R(D) \cup B^{\partial}(D) \to V^{\partial}_R(D') \cup B^{\partial}(D')$ is the induced bijection given by $g^{\partial}(c^{(j)}) = (g(c))^{(j)}$ for $c \in V_R(D)$ and     $g^{\partial}(\gamma^{(j)}) = (g(\gamma))^{(j)}$ for $\gamma \in B(D)$.
\end{enumerate}
\end{definition}

\begin{lemma}\label{lemma:Same Gauss Data}
 Two oriented twisted virtual doodle diagrams have the same Gauss data if and only if 
 they are related by a finite sequence of planar isotopies, and V1, V2, V3, V4 and T1 moves.  
\end{lemma}

\begin{figure}[ht]
  \centering
    \includegraphics[width=12cm]{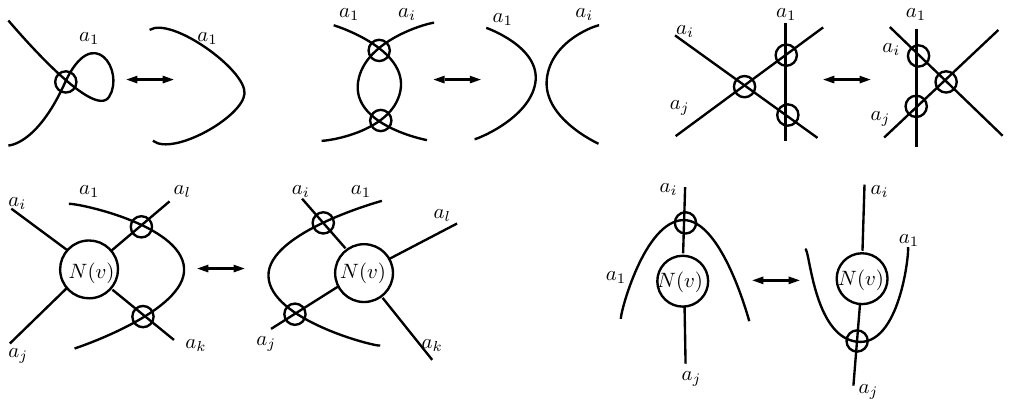}
        \caption{V1, V2, V3, V4 and T1 moves.}
        \label{mv}
        \end{figure}

\begin{proof}
Let $D$ and $D'$ be two oriented twisted virtual doodle diagrams such that they are related by a finite sequence of planar isotopies, and V1, V2, V3, V4 and T1 moves. This trivially implies that  $D$ and $D'$  admit the same Gauss data.
\par

Conversely, let $D$ and $D'$ be two oriented twisted virtual doodle diagrams admitting the same Gauss data, where without loss of generality, we assume that $D$ and $D'$ do not admit any trivial component (trivial loops). This shows that, upto planar isotopies, $D$ and $D'$ can be taken to be identical on the regular neighbourhood $N(v)$ for each $v \in V_R(D) \cup B(D)$. This further implies that $W(D)=W(D')$. Also, admitting the same Gauss data implies that there is a bijection between the arcs of $D\cap W(D)$ and those of $D'\cap W(D)$ with respect to the endpoints of the arcs. Using planar isotopies, we may assume that any intersection of an arc of $D\cap W(D)$ with an arc of $D'\cap W(D)$ is a transverse double point.  Let $a_1, a_2, \ldots, a_n$ and $a'_1, a'_2, \ldots, a'_n$ be the arcs of $D \cap W(D)$ and $D'\cap W(D)$, respectively, such that each $a_i$ maps to $a_i'$ under the above bijection. Note that each arc $a_i$ is homotopic to $a_i'$ relative to their common end points. We consider a homotopy of $a_1$ to $a_1'$ relative to the arcs $a_2, \ldots, a_n$ and the regular neighbourhoods $N(v)$ for each $v$. Any such homotopy of  $a_1$ to $a_1'$ will use a finite sequence of local moves as shown in Figure~\ref{mv}. Since all the real crossings and bars of $D$ are included in the regular neighbourhoods, we can think of the crossings in Figure~\ref{mv} as virtual crossings of $D$, and hence regard the local moves as V1, V2, V3, V4, and T1 moves. Iterating this process, we can transform each $a_i$ into $a_i'$ without changing other arcs of $D\cap W(D)$ and $D'\cap W(D)$. This proves that  $D$ and $D'$ are related by a finite sequence of planar isotopies, and V1, V2, V3, V4 and T1 moves.
\end{proof}


\par

Let $O$ be the origin of $\mathbb{R}^2$. We identify $\mathbb{R}^2 \setminus \{O\}$ with $\mathbb{R}_+\times \mathbb{S}^1$ via polar coordinates, where $\mathbb{R}_+$ is the set of positive real numbers. Let $\pi:\mathbb{R}^2 \setminus \{O\}=\mathbb{R}_+\times \mathbb{S}^1 \to \mathbb{S}^1$ denote the radial projection.  

\begin{definition}
A {\it closed twisted virtual twin diagram} is a twisted virtual doodle diagram $D$  satisfying the following conditions:
\begin{enumerate}
    \item $D$ is contained in $\mathbb{R}^2 \setminus \{O\}$.
    \item If $k: \sqcup \, \mathbb{S}^1 \to \mathbb{R}^2 \setminus \{O\}$ is the underlying immersion of $D$, then $\pi \, k: \sqcup \, \mathbb{S}^1 \to \mathbb{S}^1$ is a finite-sheeted covering map. If $D$ is oriented, we assume that $\mathbb{S}^1$ is oriented in such a way that $\pi \, k$ is orientation-preserving.
\end{enumerate}
\end{definition}

\begin{definition}
The {\it closure} of a twisted virtual twin diagram $b$ on the plane is the twisted virtual doodle diagram obtained by joining the boundary points $\{(0, 1),\ldots, (0, n) \}$ with $\{(1, 1), \ldots, (1, n) \}$ in a one-to-one fashion through non-intersecting arcs that are disjoint from $b$.
\end{definition}

\begin{remark}\label{ctwtd and closure of twtd}
The closure of a twisted virtual twin diagram is a closed twisted virtual twin diagram. Conversely, let $D$ be a closed twisted virtual twin diagram such that $\pi \, k: \sqcup \, \, \mathbb{S}^1 \to \mathbb{S}^1$ is an $n$-sheeted covering map. Consider a point $\theta \in \mathbb{S}^1$ such that $\pi^{-1}(\theta) \cap (V_R(D) \cup B(D))$ is empty. Then cutting along the ray $\pi^{-1}(\theta)$ gives a twisted virtual twin diagram $b_D$ on $n$ strands, such that the closure of $b_D$ is equivalent to $D$.
\end{remark}

\begin{figure}[ht]
  \centering
    \includegraphics[width=7cm]{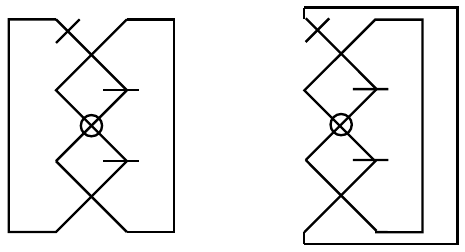}
        \caption{Two closures of the twisted virtual twin diagram representing $\gamma_1s_1\gamma_2\rho_1\gamma_2s_1$.}
        \label{one}
        \end{figure}

There are multiple ways of taking the closure of a twisted virtual twin diagram (see Figure~\ref{one}).

\begin{lemma}
Any two closures of a twisted virtual twin diagram on the plane gives equivalent twisted virtual doodle diagrams on the plane.
\end{lemma}

\begin{proof}
Let $b$ be a twisted virtual twin diagram, and let $D$ and $D'$ denote its two distinct closures. We fix an orientation on the strands of $b$ and consider $D$ and $D'$ with induced orientations. Observe that $\mu(D)=\mu(D')$, $V_R(D) = V_R(D')$ and $B(D) = B(D')$. We see that the identity map $V_R(D) \cup B(D) \to V_R(D') \cup B(D')$ satisfies the conditions for the equality of the Gauss data of $D$ and $D'$. Thus, by Lemma~\ref{lemma:Same Gauss Data}, they are related by a finite sequence of V1, V2, V3, V4 and T1 moves. Consequently, $D$ and $D'$ represent equivalent twisted virtual doodle diagrams on the plane.
\end{proof}

In view of the preceding lemma, we shall use the term the closure of a twisted virtual twin. We now prove the Alexander theorem for twisted virtual doodles.

\begin{theorem}\label{theoremAlexanderB}
Every twisted virtual doodle is equivalent to the closure of a twisted virtual twin diagram.
\end{theorem}

\begin{proof}
Let $D$ be a twisted virtual doodle diagram equipped with a fixed orientation. We construct a closed twisted virtual twin diagram $D'$ with the same Gauss data as that of $D$. We label each real crossing and bar of $D$ as in Figures \ref{r} and \ref{b}. Next, we consider the complement $\mathbb{R}^2 \setminus {\rm Int}(\mathbb{D}^2)$ of the open 2-disk in the plane and orient its boundary circle $\mathbb{S}^1$. Further, we consider the real crossings and the bars of $D$ equipped with the information of their boundary points as in Figures \ref{r} and \ref{b}. We then place them in $\mathbb{R}^2 \setminus {\rm Int}(\mathbb{D}^2)$ such that no two  elements of $V_R(D) \cup B(D)$ lie in $\pi^{-1}(\theta)$ for any $\theta \in \mathbb{S}^1$ and such that their orientations are compatible with the orientation of $\mathbb{S}^1$.  Next, we join points of $V^{\partial}_R(D) \cup B^{\partial}(D)$ by arcs in $\mathbb{R}^2\setminus {\rm Int}(\mathbb{D}^2)$ as per the Gauss data of $D$. In other words, for each $(a,b) \in X(K)$, the orientation of the arc joining $a$ to $b$ is compatible with the orientation of $\mathbb{S}^1$, and whenever two such arcs intersect, then their intersection point is marked as a virtual crossing. We see that the twisted virtual doodle diagram $D'$ obtained in this manner is a closed twisted virtual twin diagram admitting the same Gauss data as that of $D$. It now follows from  Lemma~\ref{lemma:Same Gauss Data} that $D$ and $D'$ are equivalent as twisted virtual doodle diagrams. Finally, cutting along $\pi^{-1} (\theta)$ for some $\theta \in \mathbb{S}^1$ such that $\pi^{-1} (\theta)$ does not pass through any crossing of $D'$ gives the desired twisted virtual twin diagram whose closure is $D'$, and hence equivalent to $D$. 
\end{proof}
\medskip


\section{Markov theorem for twisted virtual doodles} \label{sect:Markov}

In this section,  we prove a Markov type theorem for twisted virtual doodles. To proceed, we define various moves on twisted virtual twin diagrams.
\par

\begin{itemize}
\item A {\it twisted Markov move of type $0$} (TM0 move) is the replacement of a twisted virtual twin diagram $b$ by another diagram $b'$ such that $b$ and $b'$ are equivalent as twisted virtual twin diagrams.

\item A {\it twisted Markov move of type $1$} (TM1  move) is the replacement of a twisted virtual twin diagram $b$ by $b' b b'^{-1}$, where $b'$ is any twisted virtual twin diagram on the same number of strands. Note that, this is simply the conjugation in the twisted virtual twin group. Furthermore, the TM1  move allows replacement of a twisted virtual twin diagram of the form $b = b_1 b_2$ by $b' = b_2 b_1$. See Figure~\ref{tm1} for an illustration.
\begin{figure}[ht]
  \centering
    \includegraphics[width=10cm]{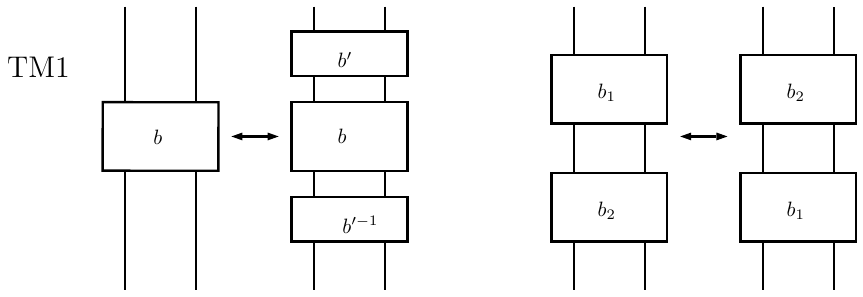}
        \caption{A twisted Markov move of type 1.}
        \label{tm1}
\end{figure}

\item For a twisted virtual twin diagram $b$ on $n$ strands and non-negative integers $s$ and $t$, let $\iota_s^t(b)$ denote the twisted virtual twin diagram on $n + s + t$ strands, which is obtained from $b$ by adding $s$ trivial strands to the left and $t$ trivial strands to the right of $b$.
This defines a group monomorphism $\iota_s^t: TVT_n \to TVT_{n+s+t}$.  
\item[] A {\it right stabilization of real or virtual type} is the replacement of a twisted virtual twin diagram $b$ on $n$ strands by $\iota_0^1(b)s_n$ or $\iota_0^1(b)\rho_n$, respectively.
\item[] Similarly, a {\it left stabilization of real or virtual type} is the replacement of a twisted virtual twin diagram $b$ on $n$ strands by $\iota_1^0(b)s_1$ or $\iota_1^0(b)\rho_1$, respectively.
\item[] A {\it twisted Markov move of type $2$} (TM2 move) is a  right or a left stabilization of real or virtual type.   See Figure~\ref{tm2}.
\begin{figure}[ht]
  \centering
    \includegraphics[width=14cm]{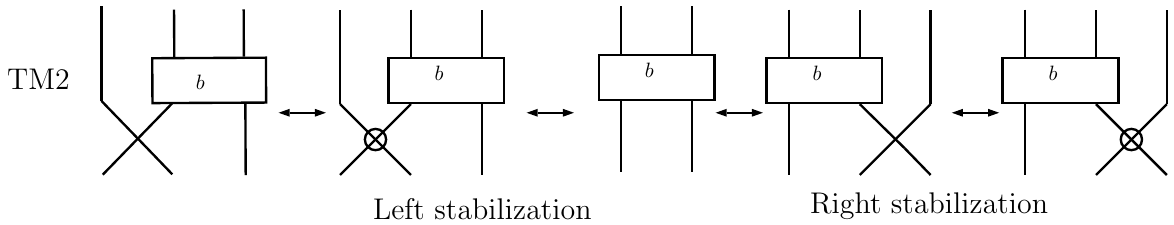}
        \caption{A twisted Markov move of type 2.}
        \label{tm2}
\end{figure}

\item  A {\it right virtual exchange move} is the replacement 
$$ \iota_0^1(b_1) s_n \iota_0^1(b_2) s_n \longleftrightarrow 
\iota_0^1(b_1) \rho_n \iota_0^1(b_2) \rho_n, $$ 
whereas a {\it left virtual exchange move} is the replacement 
$$ \iota_1^0(b_1) s_1 \iota_1^0(b_2) s_1 \longleftrightarrow   
\iota_1^0(b_1) \rho_1 \iota_1^0(b_2) \rho_1, $$ 
where $b_1$ and $b_2$ are twisted virtual twins diagrams.
\item[] A {\it twisted Markov move of type $3$} (TM3 move) is a right or a left virtual exchange move.  See Figure~\ref{tm3}.
\begin{figure}[ht]
  \centering
    \includegraphics[width=12cm]{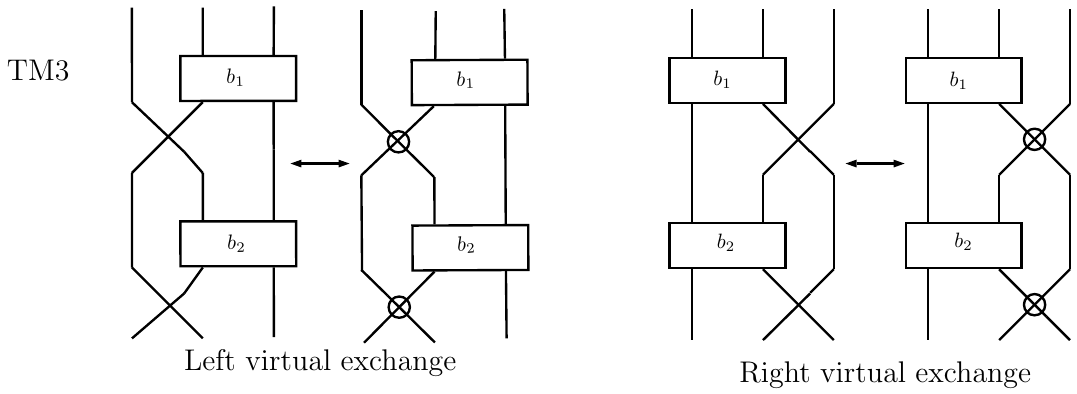}
        \caption{A twisted Markov move of type 3.}
        \label{tm3}
\end{figure}
\end{itemize}

\begin{definition}
Two twisted virtual twin diagrams are said to be {\it Markov equivalent} if they are related by a finite sequence of TM0, TM1, TM2 and TM3 moves, and planar isotopies. 
\end{definition}

\begin{lemma}\label{lem:unique}
If $D$ and $D'$ are oriented closed twisted virtual twin diagrams with the same Gauss data, then the corresponding twisted virtual twin diagrams $b_{D}$ and $b_{D'}$ are Markov equivalent. In fact, $D$ and $D'$ are related by TM0 and TM2 moves.
\end{lemma}

\begin{proof}
Let $D$ and $D'$ be oriented closed twisted virtual twin diagrams with the same Gauss data. By Lemma~\ref{lemma:Same Gauss Data}, $D$ and $D'$ are related by a finite sequence of V1, V2, V3, V4 and T1 moves. Let $\pi \, k: \sqcup \, \mathbb{S}^1 \to \mathbb{S}^1$ and $\pi \, k': \sqcup \, \mathbb{S}^1 \to \mathbb{S}^1$ be the finite-sheeted covering maps corresponding to $D$ and $D'$, respectively. For both $D$ and $D'$, using planar isotopies, we can assume that no two real crossings and/or intersection points of bars with the arcs lie on the same sheet of the covering. Let $N_1, \ldots, N_m$ and $N'_1, \ldots, N'_m$ be the regular neighbourhoods of the real crossings and the bars of $D$ and $D'$, respectively.
\par
Case 1. Suppose that $\pi(N_1), \ldots, \pi(N_m)$ and $\pi(N'_1),\ldots, \pi(N'_m)$ appear in $\mathbb{S}^1$ in the same cyclic order. Using planar isotopies, we may assume that $N_1=N'_1, \ldots, N_m=N'_m$, and that $D$ and $D'$ are identical on these regular neighbourhoods. This implies that $W(D)=W(D')$. Let $a_1, \dots, a_s$ be the arcs of $D \cap W(D)$, and $a'_1, \dots, a'_s$ be the corresponding arcs of $D '\cap W(D)$. Let $\theta \in \mathbb{S}^1$ be such that $\pi^{-1}(\theta)$ is disjoint from $N_1 \cup \dots \cup N_m$ and suppose that exists an $1\le i \le s$ such that $|a_i \cap \pi^{-1}(\theta)|\neq |a'_i \cap \pi^{-1}(\theta)|$. We move a segment of $a_i$ or $a'_i$ towards the origin by some V2 moves (which corresponds to TM0 moves of the twisted virtual twin diagram $b_D$) and apply some TM2 moves of virtual type on $b_D$ such that that $|a_i \cap \pi^{-1}(\theta)|=|a'_i \cap \pi^{-1}(\theta)|$ after these moves. Thus, we can assume that $|a_i \cap \pi^{-1}(\theta)|=|a'_i \cap \pi^{-1}(\theta)|$ for all $1\le i \le s$. Taking a homotopy from $a_1$ to $a_1'$ relative to the other arcs of $D \cap W(D)$, $D' \cap W(D)$ and the regular neighbourhoods $N_1, \ldots, N_m$, we see that $a_1$ can be transformed to $a'_1$ by a sequence of TM0  moves.  Applying this procedure inductively, we can transform the arcs $a_1, \dots, a_s$ to $a'_1, \dots, a'_s$ by a finite sequence of TM0 and TM2 moves. Thus, $b_D$ can be transformed to $b_{D'}$ by a finite sequence of TM0 and TM2 moves.
\par

Case 2. Suppose that $\pi(N_1), \dots, \pi(N_m)$ and $\pi(N'_1), \dots, \pi(N'_m)$ do not appear in $\mathbb{S}^1$ in the same cyclic order. We claim that we can interchange the positions of any two consecutive regular neighbourhoods.  Without loss of generality,  it suffices to show that positions of $\pi(N_1)$ and $\pi(N_2)$ can be interchanged. If $N_2$ is the regular neighbourhood of a real crossing, then  Figure~\ref{c} shows a procedure to interchange $\pi(N_1)$ and $\pi(N_2)$ by TM0  and TM2 moves.  
  \begin{figure}[ht]
  \centering
    \includegraphics[width=12cm]{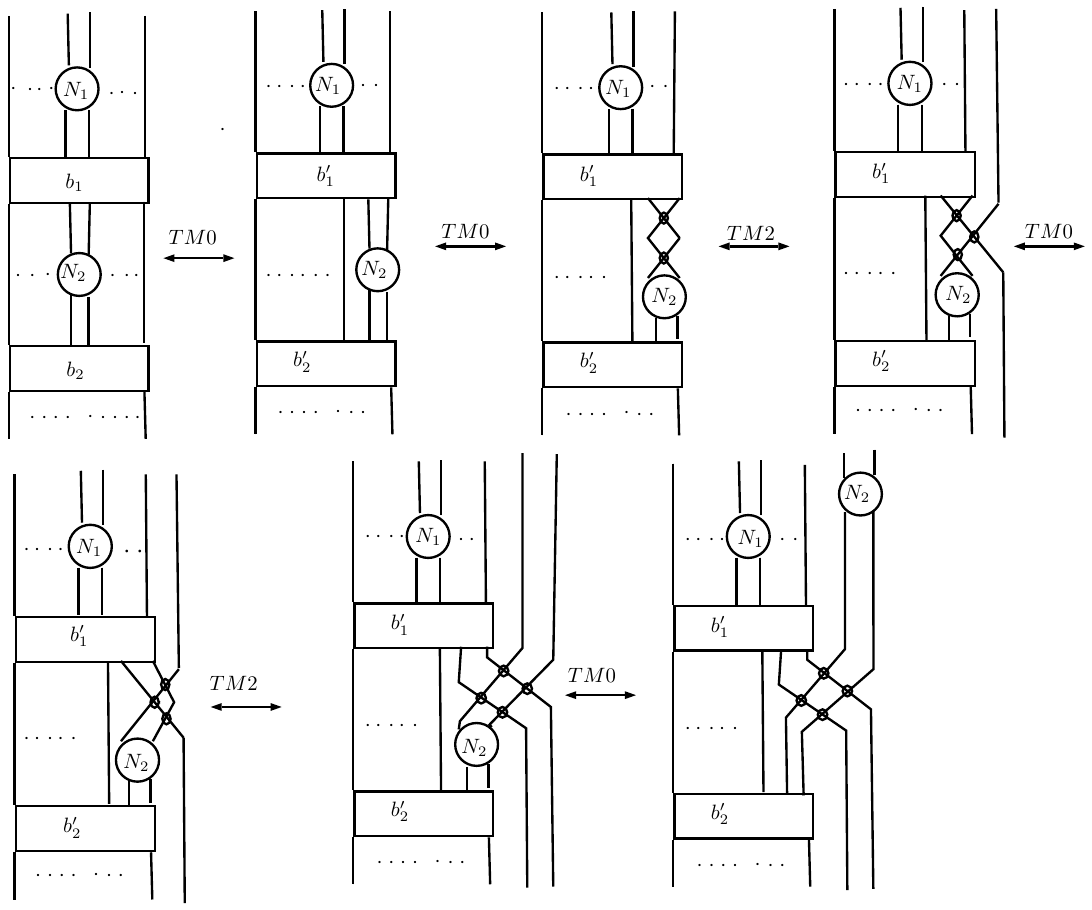}
        \caption{Interchanging the position of $N_1$ with the regular neighbourhood $N_2$ of a real crossing.}
        \label{c}
\end{figure}
If $N_2$ is the regular neighbourhood of a bar, then  Figure~\ref{c2} shows a procedure to interchange $\pi(N_1)$ and $\pi(N_2)$ by TM0  and TM2 moves.  
\begin{figure}[ht]
  \centering
    \includegraphics[width=11cm]{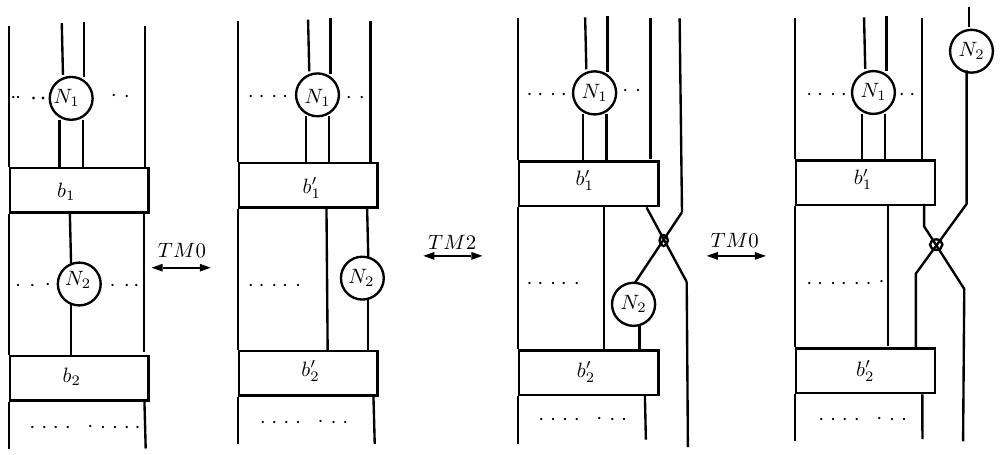}
\caption{Interchanging the position of $N_1$ with the regular neighbourhood  $N_2$ of a bar.}
        \label{c2}
\end{figure}
 This proves the claim, and the proof now follows from Case 1.
\end{proof}

We are now ready to prove the Markov theorem for twisted virtual doodles.

\begin{theorem}\label{theorem:MarkovA}
Two twisted virtual twin diagrams have equivalent closures (as twisted virtual doodles) if and only if they are Markov equivalent.
\end{theorem}

\begin{proof} 
Clearly, if two twisted virtual twin diagrams are Markov equivalent, then they have equivalent closures as twisted virtual doodles. Conversely, suppose that $D$ and $D'$ are oriented closed twisted virtual twin diagrams which are equivalent as oriented twisted virtual doodles.  Then, there is a finite sequence of oriented twisted virtual doodle diagrams $D_0=D, D_1, \ldots, D_{n-1},D_n=D'$ such that $D_{i+1}$ is obtained from $D_{i}$ by one of the oriented twisted Reidemeister moves.  
 \par
 
In general, for $1\le i \le n-1$, $D_i$ may not be a closed twisted virtual twin diagram. However, by Theorem \ref{theoremAlexanderB}, there exists a closed twisted virtual twin diagram  $\widetilde D_i$, which is equivalent to $D_i$ and has the same Gauss data as $D_i$. We assume that $D_0=\widetilde D_0$ and $D_n =\widetilde D_n$. It is sufficient to prove that $b_{\widetilde D_i}$ and $b_{\widetilde D_{i+1}}$ are Markov equivalent for each $i$.
 \par
It is proved in \cite[\textcolor{blue}{Theorem 6.6}]{MR4209535} that when $D_{i+1}$ is obtained from $D_{i}$ by a R1, R2, V1, V2, V3 or V4 move, then  $b_{\widetilde D_i}$ and $b_{\widetilde D_{i+1}}$ are Markov equivalent. Thus, it remains to consider the cases when $D_{i+1}$ is obtained from $D_{i}$ by a T1, T2 or T3 move.
\par

Case 1. Suppose that  $D_{i+1}$ is obtained from $D_i$ by a T1 move.  Then $D_{i}$ and $D_{i+1}$ have the same Gauss data, and hence $\widetilde D_{i}$ and $\widetilde D_{i+1}$ also have the same Gauss data. It follows from Lemma~\ref{lem:unique} that $b_{\widetilde D_i}$ and $b_{\widetilde D_{i+1}}$ are Markov equivalent.
\par

Case 2. Suppose that $D_{i+1}$ is obtained from $D_i$ by a T2 move. Without loss of generality, we can assume that a pair of bars is removed from $D_i$ by a T2 move to obtain $D_{i+1}$. Let $N$ be a regular neighbourhood of the arc of $D_i$ containing the two bars such that $N \cap D_i$ is an arc $\alpha$ with two bars and $N \cap D_{i+1}$ is an arc without any bars. Let $ D_i'$ be a closed twisted virtual twin diagram obtained from $D_i$ by the braiding process of Theorem \ref{theoremAlexanderB}  such that $N$ is point-wise fixed during the braiding process. Let $ D_{i+1}'$ be a closed twisted virtual twin diagram obtained from $ D_i'$ by removing the two bars from the arc $\alpha$.  Then $b_{D_i'}$ and $b_{D_{i+1}'}$ are related by a TM0  move, and hence are Markov equivalent.   Since $\widetilde D_i$ and $ D_i'$ have the same Gauss data, by Lemma~\ref{lem:unique}, $b_{\widetilde D_i}$ and $b_{D_i'}$ are Markov equivalent. Similarly, since $ D_{i+1}'$ and $\widetilde D_{i+1}$ have the same Gauss data, it follows that $b_{D_{i+1}'}$ and $b_{\widetilde D_{i+1}}$ are Markov equivalent. Thus, $b_{\widetilde D_i}$ and $b_{\widetilde D_{i+1}}$ are Markov equivalent, which is desired.
\par

Case 3.  Suppose that $D_{i+1}$ is obtained from $D_i$ by a T3 move. There are two oriented T3 moves, say, T3a and T3b as shown in Figure~\ref{ot3m}. 
\begin{figure}[ht]
  \centering
    \includegraphics[width=9cm]{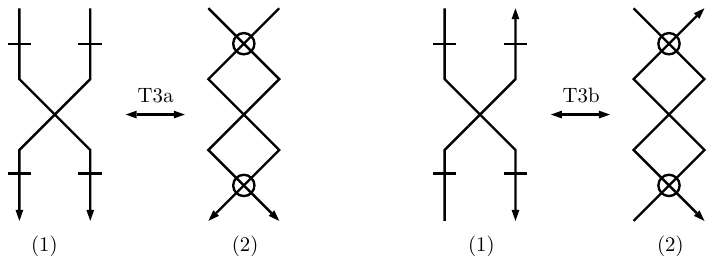}
        \caption{Oriented T3 moves.}
        \label{ot3m}
\end{figure}

First, we consider the case when $D_{i+1}$ is obtained from $D_i$ by a T3a move.  Without loss of generality, we assume that $D_i$ is the LHS of the T3a move, while $D_{i+1}$ is the RHS of the T3a move as shown in Figure~\ref{ot3m}. Let $N$ be a regular neighbourhood of the local picture of $D_i$ represented by Figure~\ref{ot3m}.  Then $N \cap D_i$ is a pair of arcs intersecting transversely at a real crossing and admitting four bars.  Let $D_i'$ be a closed twisted virtual twin diagram obtained from $D_i$ by the braiding process of Theorem \ref{theoremAlexanderB} such that $N$ is point-wise fixed during the process. Let $D_{i+1}'$ be a closed twisted virtual twin diagram obtained from $D_i'$ by applying a T3a move. Then $b_{D_i'}$ and $b_{D_{i+1}'}$ are related by a TM0  move.  Since $\widetilde D_i$ and $D_i'$ have the same Gauss data, by Lemma~\ref{lem:unique}, $b_{\widetilde D_i}$ and $b_{D_i'}$  are Markov equivalent. Similarly, since $D_{i+1}'$ and $\widetilde D_{i+1}$ have the same Gauss data,  $b_{D_{i+1}'}$ and $b_{\widetilde D_{i+1}}$ are Markov equivalent. Thus, $b_{\widetilde D_i}$ and $b_{\widetilde D_{i+1}}$ are Markov equivalent.
\par

Next, we consider the case when $D_{i+1}$ is obtained from $D_i$ by a  T3b move. We claim that a T3b move is a consequence of a T3a move together with V1, V2, V3 and V4 moves. To see this, we rotate the two diagrams of T3b move by 90 degrees clockwise. Then the LHS of T3b becomes the LHS of T3a. Further, the RHS of T3b after the rotation has a real crossing and no bars, so does the RHS of T3a.  Thus,  the two twisted virtual doodle diagrams have the same Gauss data. By Lemma \ref{lemma:Same Gauss Data}, the RHS of T3b after the rotation can be transformed to the RHS of T3a by a finite sequence of  V1, V2, V3 and V4 moves. Thus, this case reduces to that of the case of T3a move and the V1, V2, V3 and V4 moves.  
\end{proof}

\begin{remark}
By \cite[Figure 13]{MR4209535}, a left virtual stabilization move can be derived from the other Markov moves. Specifically, a left virtual stabilization move follows from a sequence of TM0, TM1 and TM2 (right virtual stabilization) moves.
\end{remark}

\begin{prop} 
If two twisted virtual twin diagrams are related by a left virtual exchange move,  then they are related by a finite sequence of TM0 and TM1 moves, and a right virtual exchange move. 
\end{prop}

\begin{proof}
Let $b$ and $b'$ be twisted virtual twin diagrams on $n$ strands that related by a left virtual exchange move. Thus, we can write
$$b= \iota_1^0(b_1) s_1 \iota_1^0(b_2) s_1 \quad \mbox{and} \quad b'= \iota_1^0(b_1) \rho_1 \iota_1^0(b_2) \rho_1$$ 
for some twisted virtual twin diagrams  $b_1$ and $b_2$ on $n-1$ strands. Let $f_n: TVT_n \to TVT_n$ be the group isomorphism defined on generators by
\begin{eqnarray*}
s_i & \mapsto& s_{n-i} \quad \text{for }  1\le i \le n-1, \\ 
\rho_i & \mapsto& \rho_{n-i} \quad \text{for }  1\le i \le n-1~\textrm{and} \\ 
\gamma_i & \mapsto& \gamma_{n-i+1} \quad \text{for }  1\le i \le n. 
\end{eqnarray*} 
Then, we have
$$f_n(b) = \iota_0^1( f_{n-1}(b_1) ) s_n \iota_0^1( f_{n-1}(b_2) ) s_n 
\quad \mbox{and} \quad 
f_n(b') =  \iota_0^1( f_{n-1}(b_1) ) \rho_n \iota_0^1( f_{n-1}(b_2) ) \rho_n, $$  and hence $f_n(b)$ and $f_n(b')$ are related by a right virtual exchange move. Let $\nabla_n$ be the twisted virtual twin diagram given by 
\begin{align*}
\nabla_n  =  \prod_{i=1}^{n-1} (\rho_i \rho_{i-1} \dots \rho_1) \prod_{j=1}^{n} \gamma_j. 
\end{align*}
\begin{figure}[ht]
  \centering
    \includegraphics[width=2cm]{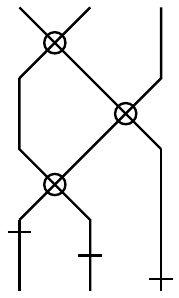}
        \caption{$\nabla_3= (\rho_1)(\rho_2\rho_1)(\gamma_1\gamma_2\gamma_3)$.}
        \label{ot3m}
\end{figure}

For instance, see Figure \ref{ot3m} for $\nabla_3$. A direct check using relations in $TVT_n$ or the diagram of $\nabla_n$ (which corresponds to TM0 moves) shows that 
$$f_n(w)=\nabla_n  w \nabla_n^{-1}$$
for all $w \in TVT_n$. In particular, the diagrams of $w$ and $f_n(w)$ are related by TM0 and TM1 moves.  Since $b$ is conjugate to $f_n(b)$ and  $b'$ is conjugate to $f_n(b')$, we conclude that $b$ and $b'$ are related by TM0  and TM1 moves, and a right virtual exchange move.
\end{proof}

\begin{remark}
In view of the preceding result, we can remove the left virtual exchange move from the definition of Markov equivalence. 
\end{remark}
\medskip

\section{Pure twisted virtual twin group}\label{sec pure}
For each $n \ge 1$, let $S_n$ be symmetric group on $n$ symbols. We consider the epimorphism  $\varphi_P: TVT_n \to S_n$ defined on generators by
$$ s_i  \mapsto  (i,i+1), \quad \rho_i  \mapsto  (i,i+1) \quad \textrm{and} \quad \gamma_j  \mapsto  1$$
for all $1\le i\le n-1$ and $1\le j\le n$. The kernel $PTVT_n$ of $\varphi_P$ is called the {\it pure twisted virtual twin group}. This gives rise to the short exact sequence
$$
1 \to PTVT_n \to TVT_n \to S_n \to 1,
$$
which splits via the section $(i, i+1) \mapsto \rho_i$, and hence $TVT_n\cong PTVT_n \rtimes S_n$.
\par

By definition, each generator $\gamma_j$ lies in $PTVT_n$. In fact, $PTVT_1= \langle \gamma_1 \rangle \cong\; \mathbb{Z}_2$. Recall from \cite[Section 5]{MR4027588} that, for each $n \ge 1$, the {\it virtual twin group}  $VT_n$ is defined as the group generated by
$$\{s_i, \rho_i \mid 1 \le i \le n-1 \}$$
subject to the following defining relations:
    \begin{align*}
        s_i^2 & = 1  & \text{ for } & 1\le i \le n-1,\\
        s_i s_j & = s_j s_i  & \text{ for } & |i-j| > 1,\\
        \rho_i^2 & = 1  & \text{ for } & 1\le i \le n-1,\\
        \rho_i \rho_j & = \rho_j \rho_i & \text{ for } & |i-j| > 1, \\
        \rho_i \rho_{i+1} \rho_i & = \rho_{i+1} \rho_i \rho_{i+1} & \text{ for } & 1 \le i \le n-2,\\
        s_i \rho_j & = \rho_j s_i &  \text{ for } & |i-j| >1,\\
        \rho_i s_{i+1} \rho_i & = \rho_{i+1} s_i \rho_{i+1} & \text{ for } & 1 \le i \le  n-2.
    \end{align*}

As above, there is an epimorphism  $VT_n \to S_n$ defined on generators by 
$$ s_i  \mapsto  (i,i+1) \quad \textrm{and} \quad \rho_i  \mapsto  (i,i+1).$$
The kernel  $PVT_n$ of this epimorphism is called the {\it pure virtual twin group}. It is known from ~\cite[Theorem 3.1]{MR4651964} that the elements
$$\lambda_{i,i+1}=s_i\rho_i,~\textrm{where}~1 \le  i\le n-1$$
and
$$\lambda_{i,j}=\rho_{j-1}\rho_{j-2}\cdots\rho_{i+1}\lambda_{i,i+1}\rho_{i+1}\cdots\rho_{j-2}\rho_{j-1},~\textrm{where}~1 \leq i < j \leq n~\text{and }~j\neq i+1$$
 generate $PVT_n$. 
\par

We use the standard presentation of $TVT_n$ and the Reidemeister--Schreier method \cite[Theorem 2.6]{MR0207802} to determine a presentation for $PTVT_n$.  For each $1 \leq k \leq n-1$, we set
$$
m_{k,i_k} :=
\begin{cases}
\rho_k \rho_{k-1} \cdots \rho_{i_k+1} &\textrm{for}~ 0 \leq i_k < k,\\[3pt]
1 &\textrm{for}~ i_k = k,
\end{cases}
$$
and
$$
M_n := 
\big\{\, m_{1,i_1} m_{2,i_2} \cdots m_{n-1,i_{n-1}} 
\;\big|\; 0 \leq i_k \leq k \text{ for each } 1 \leq k \leq n-1 \big\}
$$
as the Schreier set of coset representatives of $PTV T_n$ in $TVT_n$. For $w \in TVT_n$, let $\overline{w}$ denote the unique coset representative of the coset of $w$ in the Schreier set $M_n$.

\begin{lemma}\label{genset}
    The group $PTVT_n$ is generated by the set 
    $$\{\lambda_{k,l} \mid 1 \leq k < l \leq n \} \cup \{\gamma_j \mid 1\leq j\leq n \}$$
\end{lemma}
\begin{proof}
By \cite[Theorem 2.7]{MR0207802}, the group $PTVT_n$ is generated by the set
$$
\big\{   s_{\mu,a}:=(\mu a)  (\overline{\mu a})^{-1} \mid \mu \in M_n~\textrm{and}~a \in \{\rho_1, \ldots, \rho_{n-1}, s_1, \ldots, s_{n-1}, \gamma_1, \ldots, \gamma_ n\} \big\}.
    $$
Direct computations show that
\begin{eqnarray*}
s_{\mu, \rho_i} &= &1,\\
s_{\mu, s_i} &= & \mu (s_i\rho_i )\mu^{-1}=\mu (\lambda_{i,i+1})\mu^{-1},\\
s_{\mu,\gamma_i} &= & \mu \gamma_i \mu^{-1}.
\end{eqnarray*}

It follows from \cite[Theorem 3.1]{MR4651964} that the generators $s_{\mu, s_i}$ lies in the set $\{\lambda_{k,l} \mid 1 \leq k < l \leq n \}$. We claim that the generators $s_{\mu,\gamma_i}$ lies in the set $\{\gamma_j \mid 1\leq j\leq n \}.$  For this, we analyse the conjugation action of $S_n \cong \langle \rho_1, \ldots, \rho_{n-1} \rangle$ on the set $\{\gamma_j \mid 1\leq j\leq n \}.$  For $1\leq j, k \leq n $, we have the following cases:
\begin{itemize}
    \item If $k=j$, then
    $\rho_k\gamma_j\rho_k =\rho_j\gamma_{j}\rho_j=\gamma_{j+1}$,
    \item If $k=j-1$, then
    $\rho_k\gamma_j\rho_k =\rho_{j-1}\gamma_{j}\rho_{j-1}=\gamma_{j-1}$,
    \item If $k\neq j, j-1$, then $\rho_k\gamma_j\rho_k=\gamma_{j}$.
\end{itemize}
\par
Hence, $PTVT_n$ is generated by the set $\{\lambda_{k,l} \mid 1 \leq k < l \leq n \} \cup \{\gamma_j \mid 1\leq j\leq n \}.$
\end{proof}

Lemma~\ref{genset} and \cite[Remark 3.2]{MR4651964} leads to the following remark.

\begin{remark}\label{rem:3.2}
The action of $S_n$ on the generators  of $PTVT_n$ is given by
$$
\rho_i :
\begin{cases}
\lambda_{i,i+1} \longleftrightarrow \lambda^{-1}_{i,i+1},\\[3pt]
\lambda_{k,l} \longleftrightarrow \lambda_{k',l'} & ~\text{for}~ (k,l) \ne (i,i+1),\\[3pt]
\gamma_i \longleftrightarrow \gamma_{i+1} \\[3pt]
\gamma_{k} \longleftrightarrow \gamma_{k} & ~\text{for}~ k \neq i, i+1
\end{cases}
$$
where the transposition $(k',l')$ equals to $(i,i+1)(k,l)(i,i+1)$ and $k' < l'$.
\end{remark}

\begin{theorem}\label{sTVP_n}
For $n \ge 2$,  the group $PTVT_n$ admits a presentation with generating set 
    $$\{\lambda_{k,l} \mid 1 \leq k < l \leq n \} \cup \{\gamma_j \mid 1\leq j\leq n \}$$ and defining relations given by
     \begin{align}
\lambda_{i,j}\lambda_{k,l} &=\lambda_{k,l}\lambda_{i,j},\label{comm-clas}\\
\gamma_i^2 &=1,\label{g1}\\
\gamma_i\gamma_j &=\gamma_j\gamma_i,\label{g2}\\
\lambda_{i,j}\gamma_k & =\gamma_k\lambda_{i,j},\label{g3}\\
\lambda_{i,j}^{-1} &=\gamma_i\gamma_j\lambda_{i,j}\gamma_j\gamma_i,\label{g4}
\end{align}
where distinct letters stand for distinct indices.
\end{theorem}

\begin{proof}
In view of Lemma~\ref{genset}, it remains to determine the defining relations. Given an element $u \in PTVT_n$ which is written in terms of the  generators of $TVT_n$, we define a rewriting process $\tau$ that writes $u$ in  the generators of $PTVT_n$. To a word 
$$
u=a_1^{\epsilon_1}a_2^{\epsilon_2}\cdots a_v^{\epsilon_v},
$$
where $\epsilon_l=\pm 1$ and $a_l \in \{s_1,s_2, \ldots, s_{n-1}, \rho_1,\rho_2, \ldots, \rho_{n-1}, \gamma_1,\gamma_2, \ldots, \gamma_n\}$, consider the word
   $$\tau(u)=s_{k_1,a_1}^{\epsilon_1}s_{k_2,a_2}^{\epsilon_2}\cdots s_{k_v,a_v}^{\epsilon_v}$$
   in the generators of $PTVT_n$, where $k_j$ is a representative in ${\rm M}_n$ of the $(j-1)$-th initial segment of the word $u$ if $\epsilon_j=1$, and it is a representative of the $j$-th initial segment of $u$ if $\epsilon_j=-1$.   By \cite[Theorem 2.9]{MR0207802}, the group $PTVT_n$ has defining relations
   $$ \{\tau(\mu r \mu^{-1}) \mid \mu \in M_n~\textrm{and}~r~ \textrm{is a defining relation of}~TVT_n\}.$$
\par

Observe that 
$$\overline{\alpha} = \alpha,$$
$$\overline{\alpha_1 s_{i_1} \cdots \alpha_k s_{i_k} }
= \alpha_1 \rho_{i_1} \cdots \alpha_k \rho_{i_k},$$
$$\overline{\alpha_1 \gamma_{i_1} \cdots \alpha_k \gamma_{i_k} }
= \alpha_1 \cdots \alpha_k$$
for elements $\alpha, \alpha_j$ in the subgroup  $\langle \rho_1, \ldots, \rho_{n-1} \rangle$.  
\par

The common relations \eqref{rel-height-s2}-\eqref{rel-vsv} in the presentations of $VT_n$ and $TVT_n$ give rise to relations \eqref{comm-clas} in $PTVT_n$, which are derived in \cite[Theorem 3.3]{MR4651964}. It remains to consider the additional relations \eqref{rel-inverse-b}-\eqref{rel-twist-III} in the presentation of $TVT_n$.
Let us take $\mu = \rho_{i_1}\rho_{i_2}\cdots\rho_{i_k} \in M_n$.

\begin{itemize}
\item First consider the relation $\gamma_i^2=1$. Then we have
   \begin{align*}
\tau(\mu\gamma_i^2\mu^{-1}) & = \tau(\rho_{i_1}\rho_{i_2}\cdots\rho_{i_k}\gamma_i\gamma_i\rho_{i_k}\rho_{i_{k-1}}\cdots\rho_{i_1}) \\
& = s_{1,\rho_{i_1}}s_{\bar{\rho_{i_1}},\rho_{i_2}}\ldots s_{\bar{\mu},\gamma_i}s_{\overline{\mu\gamma_i}, \gamma_i}\ldots s_{\overline{\mu\gamma_i\gamma_i\cdots \rho_{i_k-1}}, \rho_{i_k}}\\
    & = s_{\mu,\gamma_i}s_{\mu, \gamma_i}\\
 & = (\mu\gamma_i\mu^{-1})(\mu\gamma_i\mu^{-1})\\
  & = \mu\gamma_i^2\mu^{-1}.
   \end{align*}
For $\mu=1$, we obtain $\gamma_i^2=1$. In general, by Remark~\ref{rem:3.2}, $\mu\gamma_i^2\mu^{-1}$ gives the relations \eqref{g1} in $PTVT_n$.

\item  For the relation $\gamma_i\gamma_j\gamma_i\gamma_{j}=1$, we have
 \begin{align*}
\tau(\mu\gamma_i\gamma_j\gamma_i\gamma_{j}\mu^{-1})& = s_{1,\rho_{i_1}}s_{\bar{\rho_{i_1}},\rho_{i_2}}\ldots s_{\bar{\mu},\gamma_i}s_{\overline{\mu\gamma_i},\gamma_j}s_{\overline{\mu\gamma_i\gamma_j},\gamma_i}s_{\overline{\mu\gamma_i\gamma_j\gamma_i}, \gamma_j}\ldots s_{\overline{\mu\gamma_i\gamma_j\gamma_i\gamma_j\cdots \rho_{i_k-1}}, \rho_{i_k}}\\
& = s_{\mu,\gamma_i}s_{\mu,\gamma_j}s_{\mu, \gamma_i}s_{\mu,\gamma_j}\\
& = (\mu\gamma_i\mu^{-1}) (\mu\gamma_j \mu^{-1}) (\mu\gamma_i \mu^{-1}) (\mu\gamma_j\mu^{-1})\\
& = \mu\gamma_i\gamma_j\gamma_i\gamma_j\mu^{-1}.
            \end{align*}
For $\mu=1$, we get $\gamma_i\gamma_j\gamma_i\gamma_j=1$.  By Remark~\ref{rem:3.2}, the remaining relations $\mu\gamma_i\gamma_j\gamma_i\gamma_j\mu^{-1}$ gives the relations \eqref{g2} in $PTVT_n$.

\item Next, we consider the relation $s_i\gamma_ks_i\gamma_k=1$ for $k\neq i,i+1$. Then we have
\begin{align*}
 \tau(\mu s_i\gamma_ks_i\gamma_k \mu^{-1})& = s_{1,\rho_{i_1}}s_{\bar{\rho_{i_1}},\rho_{i_2}}\ldots s_{\bar{\mu},s_i}s_{\overline{\mu s_i},\gamma_k}s_{\overline{\mu s_i\gamma_k},s_i}s_{\overline{ \mu s_i\gamma_ks_i}, \gamma_k}\ldots s_{\overline{\mu s_i\gamma_ks_i\gamma_k\cdots \rho_{i_k-1}}, \rho_{i_k}} \\
    & = s_{\mu,s_i}s_{\mu\rho_i,\gamma_k}s_{\mu\rho_i, s_i}s_{\mu\rho_i\rho_i,\gamma_k}\\
  & = (\mu\lambda_{i,i+1}\mu^{-1})(\mu\rho_i\gamma_k\rho_i\mu^{-1})(\mu\rho_i\lambda_{i,i+1}\rho_i\mu^{-1})(\mu\gamma_k\mu^{-1})\\
  & = (\mu\lambda_{i,i+1}\mu^{-1})(\mu\gamma_k\mu^{-1})(\mu\lambda^{-1}_{i,i+1}\mu^{-1})(\mu\gamma_k\mu^{-1})\\
  & = \mu(\lambda_{i,i+1}\gamma_k\lambda^{-1}_{i,i+1}\gamma_k)\mu^{-1}.
            \end{align*}
For $\mu=1$, we get $\lambda_{i,i+1}\gamma_k\lambda^{-1}_{i,i+1}\gamma_k=1$. By Remark~\ref{rem:3.2}, the remaining relations $\mu\lambda_{i,i+1}\gamma_k\lambda^{-1}_{i,i+1}\gamma_k\mu^{-1}$ gives the relations \eqref{g3} in $PTVT_n$.
\item  For the relation $\gamma_j\rho_i\gamma_{j}\rho_i=1$ with $j\neq i,i+1$, we have
 \begin{align*}
\tau(\mu\gamma_j\rho_i\gamma_{j}\rho_i\mu^{-1}) & = s_{1,\rho_{i_1}}s_{\bar{\rho_{i_1}},\rho_{i_2}}\ldots s_{\bar{\mu},\gamma_j}s_{\overline{\mu\gamma_j},\rho_i}s_{\overline{\mu\gamma_j\rho_i},\gamma_j}s_{\overline{\mu\gamma_j\rho_i\gamma_j}, \rho_i}\ldots s_{\overline{\mu\gamma_j\rho_i\gamma_j\rho_i\cdots \rho_{i_k-1}}, \rho_{i_k}}\\
& = s_{\mu,\gamma_j}s_{\mu\rho_i,\gamma_j}\\
& = (\mu\gamma_j\mu^{-1}) (\mu\rho_i\gamma_j\rho_i \mu^{-1})\\
& = (\mu\gamma_j\mu^{-1}) (\mu\gamma_j\mu^{-1})\\
& = 1.
\end{align*}
Thus, for each $\mu$, the relation $\gamma_j\rho_i\gamma_j\rho_i=1$ for $j\neq i,i+1$ gives the trivial relation.

\item  For the relation $\gamma_{i+1}\rho_i\gamma_{i}\rho_i=1$, we have
 \begin{eqnarray*}
&& 
\tau(\mu\gamma_{i+1}\rho_i\gamma_{i}\rho_i\mu^{-1})\\
& = & s_{1,\rho_{i_1}}s_{\bar{\rho_{i_1}},\rho_{i_2}}\ldots s_{\bar{\mu},\gamma_{i+1}}s_{\overline{\mu\gamma_{i+1}},\rho_i}s_{\overline{\mu\gamma_{i+1}\rho_i},\gamma_{i}}s_{\overline{\mu\gamma_{i+1}\rho_i\gamma_i}, \rho_i}\ldots s_{\overline{\mu\gamma_{i+1}\rho_i\gamma_i\rho_i\cdots \rho_{i_k-1}}, \rho_{i_k}}\\
& = & s_{\mu,\gamma_{i+1}}s_{\mu\rho_i,\gamma_i}\\
& =& (\mu\gamma_{i+1}\mu^{-1}) (\mu\gamma_{i+1} \mu^{-1})\\
& =& 1.
\end{eqnarray*}
\item For the last relation $\rho_is_i\rho_i\gamma_{i+1}\gamma_i s_i\gamma_i\gamma_{i+1}=1$, we have
\begin{eqnarray*}
&& \tau(\mu\rho_is_i\rho_i\gamma_{i+1}\gamma_i s_i\gamma_i\gamma_{i+1}\mu^{-1})\\
&= & s_{1,\rho_{i_1}}s_{\bar{\rho_{i_1}},\rho_{i_2}}\ldots s_{\bar{\mu},\rho_i}s_{\overline{\mu\rho_i},s_i}s_{\overline{\mu\rho_is_i},\rho_i},s_{\overline{\mu\rho_is_i\rho_i},\gamma_{i+1}}s_{\overline{\mu\rho_is_i\rho_i\gamma_{i+1}},\gamma_i}s_{\overline{\mu\rho_is_i\rho_i\gamma_{i+1}\gamma_i},s_i}\\
 & & s_{\overline{\mu\rho_is_i\rho_i\gamma_{i+1}\gamma_i s_i},\gamma_i}s_{\overline{\mu\rho_is_i\rho_i\gamma_{i+1}\gamma_i s_i\gamma_i},\gamma_{i+1}}\ldots s_{\overline{\mu\rho_is_i\rho_i\gamma_{i+1}\gamma_i s_i\gamma_i\gamma_{i+1}\cdots \rho_{i_k-1}}, \rho_{i_k}}\\
   &= &
s_{\mu\rho_i,s_i}s_{\mu\rho_i,\gamma_{i+1}}s_{\mu\rho_i,\gamma_{i}}s_{\mu\rho_i, s_i}s_{\mu,\gamma_i}s_{\mu,\gamma_{i+1}}\\
&= & (\mu\rho_i\lambda_{i,i+1}\rho_i\mu^{-1}) (\mu\rho_i\gamma_{i+1}\rho_i\mu^{-1})(\mu\rho_i\gamma_{i}\rho_i\mu^{-1})(\mu\rho_i\lambda_{i,i+1}\rho_i\mu^{-1})(\mu\gamma_i\mu^{-1})(\mu\gamma_{i+1}\mu^{-1})\\
&= & (\mu\lambda^{-1}_{i,i+1}\mu^{-1}) (\mu\gamma_{i}\mu^{-1})(\mu\gamma_{i+1}\mu^{-1})(\mu\lambda^{-1}_{i,i+1}\mu^{-1})(\mu\gamma_i\mu^{-1})(\mu\gamma_{i+1}\mu^{-1})\\
&= & \mu(\lambda^{-1}_{i, i+1} \gamma_{i}\gamma_{i+1}\lambda^{-1}_{i,i+1}\gamma_i\gamma_{i+1})\mu^{-1}.
            \end{eqnarray*}
For $\mu=1$, we get $\lambda^{-1}_{i+1,i} \gamma_{i}\gamma_{i+1}\lambda^{-1}_{i,i+1}\gamma_i\gamma_{i+1}=1$. Again, by Remark~\ref{rem:3.2}, the remaining relations $\mu\lambda^{-1}_{i+1,i} \gamma_{i}\gamma_{i+1}\lambda^{-1}_{i,i+1}\gamma_i\gamma_{i+1}\mu^{-1}$ gives the relations \eqref{g4} in $PTVT_n$. 
\end{itemize}
Thus, the group $PTVT_n$ is defined by the relations \eqref{comm-clas}-\eqref{g4}, and the proof is complete.
\end{proof} 
\medskip

\section{Decompositions of $TVT_n$ and $PTVT_n$}\label{sec decompositions}
In this section, we derive decompositions of $TVT_n$ and $PTVT_n$, and give some applications.
 
\begin{theorem}\label{decom1}
For each $n \geq 2$, the pure twisted virtual twin group $PTVT_n$ admits the decomposition
$$
PTVT_n \;\cong\; \big(F_\infty^{(n)} * \mathbb{Z}_2\big) \rtimes PTVT_{n-1},
$$
where $F_\infty^{(n)}$ is a free group of countably infinite rank generated by conjugates of the elements $\lambda_{i, n}$  for $1 \leq i \leq n-1$, and the free factor $\mathbb{Z}_2$ is generated by $\gamma_n$. Iterating the decomposition yields
$$
PTVT_n \;\cong\; 
\Big(F_\infty^{(n)} * \mathbb{Z}_2\Big) \rtimes
\Big(F_\infty^{(n-1)} * \mathbb{Z}_2\Big) \rtimes \cdots \rtimes
\big(\mathbb{Z} * \mathbb{Z}_2\big) \rtimes \mathbb{Z}_2.
$$
\end{theorem}

\begin{proof}
Consider the epimorphism
$$
\pi_n : PTVT_n \longrightarrow PTVT_{n-1}
$$
defined by
$$
\pi_n(\lambda_{i,j}) = 
\begin{cases}
\lambda_{i,j} &~\textrm{for}~ 1 \le i <j \leq n-1, \\
1 & ~\textrm{for}~j=n,
\end{cases}
\quad \textrm{and} \quad
\pi_n(\gamma_k) =
\begin{cases}
\gamma_k &~\textrm{for}~ 1 \le k \leq n-1, \\
1 &~\textrm{for}~ k=n.
\end{cases}
$$
Geometrically, the map $\pi_n$ simply forgets the $n$-th strand in the twisted virtual twin diagram representing an element of $PTVT_n$. Let
$$
1 \;\to\; \ker(\pi_n) \;\to\; PTVT_n \;\xrightarrow{\pi_n}\; PTVT_{n-1} \;\to\; 1
$$
be the associated short exact sequence. Note that, the inclusion map $\iota_n : PTVT_{n-1} \hookrightarrow PTVT_n$ satisfies $\pi_n \, \iota_n = \mathrm{id}_{PTVT_{n-1}}$, and hence
$$
PTVT_n \;\cong\; \ker(\pi_n) \rtimes PTVT_{n-1}.
$$
We use the Reidemeister--Schreier method to find a presentation of $\ker(\pi_n)$. In view of the above decomposition, we take the Schreier system to be $PTVT_{n-1}$.  For $\mu \in PTVT_{n-1}$, we have
$$
s_{\mu, \lambda_{i,j}} =
\begin{cases}
1 &~\textrm{for}~ 1 \le i <j \leq n-1,\\
\mu \lambda_{i,j} \mu^{-1} & ~\textrm{for}~j=n,
\end{cases}
\quad \textrm{and} \quad
s_{\mu,\gamma_{k} } =
\begin{cases}
1 &~\textrm{for}~ 1 \le k \leq n-1, \\
\mu \gamma_{k} \mu^{-1} &~\textrm{for}~ k=n.
\end{cases}
$$
This implies that $\ker (\pi_n)$ is generated by the set
$$
\{\mu \lambda_{i,n} \mu^{-1},~ \mu \gamma_{n} \mu^{-1}  \mid \mu \in PTVT_{n-1}~\textrm{and}~ 1 \le i  \le n-1 \}.
$$
For each $n \ge 2$, a direct check shows that the rewriting process $\tau(\mu r \mu^{-1})$ gives a trivial relation of $\ker(\pi_n)$ whenever $\mu \in PTVT_{n-1}$ and $r$ is a relation of $PTVT_n$ consisting of generators from the set  $\{\lambda_{k,l} \mid 1 \leq k < l \leq n-1 \} \cup \{\gamma_j \mid 1\leq j\leq n-1 \}$. Thus, the only possible non-trivial relations of  $\ker(\pi_n)$ arise from the relations of $PTVT_n$ that consist of at least one generator from the set  $\{\lambda_{k,n} \mid 1 \leq k \le n-1\} \cup \{\gamma_n \}$.
\par

\begin{itemize}
\item Applying the rewriting process to the relation $\gamma_n^2=1$ in $PTVT_n$ gives
\begin{eqnarray*}
\tau(\mu \gamma_n^2\mu^{-1})
&=& s_{\overline{\mu}, \gamma_n} s_{\overline{\mu \gamma_n}, \gamma_n}\\
&=& s_{\mu, \gamma_n} s_{\mu, \gamma_n}\\
&=& \mu \gamma_n^2 \mu^{-1},
\end{eqnarray*}
which yields the relation $\gamma_n^2=1$ in $\ker(\pi_n)$.

\item Applying the rewriting process to the relation  $\lambda_{i,n}\gamma_i\gamma_n\lambda_{i,n}\gamma_n\gamma_i = 1$ in $PTVT_n$ gives
\begin{eqnarray*}
\tau(\mu \lambda_{i,n}\gamma_i\gamma_n\lambda_{i,n}\gamma_n\gamma_i \mu^{-1})
&=& s_{\overline{\mu}, \lambda_{i,n}} s_{\overline{\mu\lambda_{i,n}}, \gamma_{i}}
   s_{\overline{\mu\lambda_{i,n}\gamma_{i}}, \gamma_{n}} s_{\overline{\mu\lambda_{i,n}\gamma_{i}\gamma_n}, \lambda_{i,n}}
   s_{\overline{\mu\lambda_{i,n}\gamma_{i}\gamma_n\lambda_{i,n}}, \gamma_{n}}
   s_{\overline{\mu\lambda_{i,n}\gamma_{i}\gamma_n\lambda_{i,n}\gamma_n}, \gamma_{i}}\\
&=& s_{\mu, \lambda_{i,n}} s_{\mu\gamma_{i}, \gamma_{n}}s_{\mu\gamma_i, \lambda_{i,n}} s_{\mu\gamma_{i}, \gamma_{n}} \\
&=& (\mu \lambda_{i,n}  \mu^{-1}) (\mu \gamma_{i} \gamma_{n} \gamma_i \mu^{-1}) (\mu\gamma_i \lambda_{i,n} \gamma_i \mu^{-1}) (\mu \gamma_{i}\gamma_{n}\gamma_i \mu^{-1}) \\
&=& \mu \lambda_{i,n}   \gamma_{i} \gamma_{n}  \lambda_{i,n} \gamma_{n}\gamma_i \mu^{-1}.
\end{eqnarray*}
This gives the relation
$$
 \gamma_{n}  \lambda_{i,n}^{-1} \gamma_n=  \gamma_i  \lambda_{i,n} \gamma_i
$$
in $\ker(\pi_n)$, which expresses the generator $\gamma_i  \lambda_{i,n} \gamma_i$ in terms of the generators $\gamma_{n}$ and $ \lambda_{i,n}$.

\item Next, applying the rewriting process to the relation $\gamma_i\gamma_n \gamma_i\gamma_n=1$ in $PTVT_n$ gives 
\begin{eqnarray*}
\tau(\mu \gamma_i\gamma_n \gamma_i\gamma_n \mu^{-1})
&=& s_{\overline{\mu}, \gamma_i} s_{\overline{\mu\gamma_i}, \gamma_n}
   s_{\overline{\mu \gamma_i \gamma_n}, \gamma_i} s_{\overline{\mu\gamma_i \gamma_n \gamma_i}, \gamma_n}\\
&=&  s_{\mu\gamma_{i}, \gamma_{n}}  s_{\mu, \gamma_{n}} \\
&=& \mu \gamma_{i} \gamma_n \gamma_{i} \gamma_{n}  \mu^{-1}.
\end{eqnarray*}
This yields the relation
$$
\gamma_{i} \gamma_n \gamma_{i} =\gamma_{n}
$$
in $\ker(\pi_n)$, which simply identifies two generators of $\ker(\pi_n)$.
\par

Note that $PTVT_1 = \langle \gamma_1 \rangle \cong \mathbb{Z}_2$. It follows from the preceding computations that $\ker(\pi_2)$ is generated by the set 
$$
\{\, \lambda_{1,2}, \gamma_1 \lambda_{1,2}\gamma_1,\  \gamma_{2} , \gamma_1\gamma_2\gamma_1\,\}.
$$
and has defining relations $\gamma_2^2=1$, $\gamma_1\gamma_2\gamma_1=\gamma_2$ and $\gamma_2\lambda_{1,2}^{-1}\gamma_2=\gamma_1\lambda_{1,2}\gamma_1$. Hence,  we have $$\ker(\pi_2) = \langle\lambda_{1,2}, \gamma_2~ | ~ \gamma_2^2=1 \rangle \cong \mathbb{Z} * \mathbb{Z}_2,$$ and therefore $PTVT_2 \cong (\mathbb{Z} * \mathbb{Z}_2) \rtimes \mathbb{Z}_2$.
\par

\item For $n \ge 3$, we have the additional relation $\lambda_{i,j}\gamma_k \lambda_{i,j}^{-1} \gamma_k=1$ in $PTVT_n$, where $i, j, k$ are distinct indices. Applying the rewriting process to the relation $\lambda_{i,n}\gamma_k\lambda_{i,n}^{-1}\gamma_k = 1$ in $PTVT_n$ gives
\begin{eqnarray*}
\tau(\mu \lambda_{i,n}\gamma_k\lambda_{i,n}^{-1}\gamma_k\mu^{-1})
&=& s_{\overline{\mu}, \lambda_{i,n}} s_{\overline{\mu\lambda_{i,n}}, \gamma_{k}}
   s_{\overline{\mu\lambda_{i,n}\gamma_{k} \lambda_{i,n}^{-1}}, \lambda_{i,n}}^{-1}
   s_{\overline{\mu\lambda_{i,n}\gamma_{k}\lambda_{i,n}^{-1}}, \gamma_{k}} \\
&=& s_{\mu, \lambda_{i,n}} s_{\mu\gamma_{k}, \lambda_{i,n}}^{-1} \\
&=& ( \mu \lambda_{i,n} \mu^{-1})  ( \mu \gamma_{k} \lambda_{i,n} \gamma_{k}\mu^{-1}\big)^{-1}.
\end{eqnarray*}
This gives the relation
 $$
 \lambda_{i,n}  = \gamma_{k}\lambda_{i,n}\gamma_{k}
$$
in $\ker(\pi_n)$, which simply identifies two of its generators. Similarly, applying the rewriting process to the relation $\lambda_{i,j}\gamma_n\lambda_{i,j}^{-1}\gamma_n=1$ gives
\begin{eqnarray*}
\tau(\mu \lambda_{i,j}\gamma_n\lambda_{i,j}^{-1}\gamma_n \mu^{-1})
&=& s_{\overline{\mu}, \lambda_{i,j}} s_{\overline{\mu\lambda_{i,j}}, \gamma_n}
   s_{\overline{\mu\lambda_{i,j}\gamma_n \lambda_{i,j}^{-1}}, \lambda_{i,j}}^{-1}
   s_{\overline{\mu\lambda_{i,j}\gamma_n\lambda_{i,j}^{-1}}, \gamma_n} \\
&=& s_{\mu \lambda_{i,j}, \gamma_n} s_{\mu, \gamma_n} \\
&=& ( \mu \lambda_{i,j} \gamma_n \lambda_{i,j}^{-1} \mu^{-1} ) ( \mu \gamma_n \mu^{-1}).
\end{eqnarray*}
This gives the relation
$$
 \lambda_{i,j}\gamma_n \lambda_{i,j}^{-1} = \gamma_n,
$$
which again identifies two generators of $\ker(\pi_n)$.
\par

\item For $n \ge 4$, we have the additional relation $\lambda_{i,j}\lambda_{k,l}\lambda_{i,j}^{-1}\lambda_{k,l}^{-1} = 1$ in $PTVT_n$, where $i,j,k,l$ are distinct indices. We only need to consider the case when precisely one of $i,j,k,l$ is equal to $n$. Without loss of generality, assume that $j = n$.  Then applying the rewriting process to the relation $\lambda_{i,n}\lambda_{k,l}\lambda_{i,n}^{-1}\lambda_{k,l}^{-1} = 1$ in $PTVT_n$ gives
\begin{eqnarray*}
\tau(\mu \lambda_{i,n}\lambda_{k,l}\lambda_{i,n}^{-1}\lambda_{k,l}^{-1}\mu^{-1})
&=& s_{\overline{\mu}, \lambda_{i,n} }s_{\overline{\mu\lambda_{i,n}}, \lambda_{k,l}}
   s_{\overline{\mu\lambda_{i,n}\lambda_{k,l} \lambda_{i,n}^{-1}}, \lambda_{i,n}}^{-1}
   s_{\overline{\mu\lambda_{i,n}\lambda_{k,l}\lambda_{i,n}^{-1} \lambda_{k,l}^{-1}}, \lambda_{k,l}}^{-1} \\
&=& s_{\mu, \lambda_{i,n}}s_{\mu\lambda_{k,l}, \lambda_{i,n}}^{-1} \\
&=& (\mu \lambda_{i,n} \mu^{-1}) (\mu \lambda_{k,l} \lambda_{i,n} \lambda_{k,l}^{-1} \mu^{-1})^{-1}.
\end{eqnarray*}
This gives the relation
$$
 \lambda_{i,n}  =  \lambda_{k,l}\lambda_{i,n}\lambda_{k,l}^{-1}
 $$
in $\ker (\pi_n)$, which simply identifies two of its generators.
\end{itemize}

Finally, to show that there are still infinitely many distinct generators of  $\ker (\pi_n)$ for $n \ge 3$, consider the subset
$$
\{\lambda_{n-2,n-1}^{\varepsilon}\, \lambda_{n-1,n}\, \lambda_{n-2,n-1}^{-\varepsilon} \mid \varepsilon \in \mathbb{Z} \}
$$
of generators of $\ker (\pi_n)$. It follows from the preceding computations that 
$$\lambda_{n-2,n-1}^{\varepsilon}\, \lambda_{n-1,n}\, \lambda_{n-2,n-1}^{-\varepsilon} \ne \lambda_{n-2,n-1}^{\varepsilon'}\, \lambda_{n-1,n}\, \lambda_{n-2,n-1}^{-\varepsilon'}$$ 
for each $\varepsilon \ne \varepsilon'$, and hence it is a countably infinite set. Consequently, $\ker (\pi_n)$ is generated by a countably infinite set of elements with only the element $\gamma_n$ being of order two and with no other defining relations. This implies that $\ker (\pi_n) \cong F^{(n)}_\infty* \mathbb{Z}_2$, where $F^{(n)}_\infty$ is a free group of countably infinite rank. Now, iterating the decomposition, we obtain

$$
PTVT_n \;\cong\; \Big(F_\infty^{(n)} * \mathbb{Z}_2\Big) \rtimes
\Big(F_\infty^{(n-1)} * \mathbb{Z}_2\Big) \rtimes \cdots \rtimes \Big(F_\infty^{(3)} * \mathbb{Z}_2\Big) \rtimes 
\big(\mathbb{Z} * \mathbb{Z}_2\big) \rtimes \mathbb{Z}_2,
$$
which completes the proof.
\end{proof}

Next, we consider another decomposition of $PTVT_n$. Consider the epimorphism 
$$
\psi_P \colon PTVT_n \to \mathbb{Z}_2^n
$$ 
given by $\lambda_{k,l} \mapsto 1$ and $\gamma_j \mapsto e_j$ for all $j, k, l$. If $H_n$ denotes $\ker(\psi_P)$ and $\mathbb{Z}_2^n= \langle e_1, \ldots, e_n\rangle$, then the short exact sequence
$$
1 \longrightarrow H_n \longrightarrow PTVT_n \longrightarrow \mathbb{Z}_2^n \longrightarrow 1
$$
splits via the section $e_i \mapsto \gamma_i$, and hence 
$$
PTVT_n \cong H_n \rtimes \mathbb{Z}_2^n.
$$

Let us set $a^b=bab^{-1}$.

\begin{theorem} \label{decom2}
For each $n \geq 2$, the group $H_n$ is an irreducible right-angled Artin group  generated by the set
$$\{\lambda_{i,j}, ~\lambda_{i,j}^{\gamma_i},~ \lambda_{i,j}^{\gamma_j} \mid 1\leq i<j \leq n \}$$
and admitting the following defining relations:
$$\lambda_{i,j}\lambda_{k,l}=\lambda_{k,l}\lambda_{i,j},~~ \lambda^{\gamma_i}_{i,j}\lambda_{k,l} =\lambda_{k,l}\lambda^{\gamma_i}_{i,j}, ~~\lambda^{\gamma_j}_{i,j}\lambda_{k,l} =\lambda_{k,l}\lambda^{\gamma_j}_{i,j},$$
$$\lambda_{i,j}\lambda^{\gamma_k}_{k,l} =\lambda^{\gamma_k}_{k,l}\lambda_{i,j},~~\lambda_{i,j}\lambda^{\gamma_l}_{k,l} =\lambda_{k,l}^{\gamma_l}\lambda_{i,j},~~\lambda^{\gamma_i}_{i,j}\lambda^{\gamma_k}_{k,l} =\lambda^{\gamma_k}_{k,l}\lambda^{\gamma_i}_{i,j}, $$
$$\lambda^{\gamma_i}_{i,j}\lambda^{\gamma_l}_{k,l} =\lambda^{\gamma_l}_{k,l}\lambda^{\gamma_i}_{i,j},~~ \lambda^{\gamma_j}_{i,j}\lambda^{\gamma_k}_{k,l} =\lambda^{\gamma_k}_{k,l}\lambda^{\gamma_j}_{i,j},~~\lambda^{\gamma_j}_{i,j}\lambda^{\gamma_l}_{k,l} =\lambda^{\gamma_l}_{k,l}\lambda^{\gamma_j}_{i,j},$$ 
where distinct letters stand for distinct indices.
\end{theorem}

\begin{proof}
We use the Reidemeister--Schreier method. Since $PTVT_n \cong H_n \rtimes \mathbb{Z}_2^n$,  we take the Schreier system to be $\mathbb{Z}^n_2$. 
For each $\mu \in \mathbb{Z}^n_2$ and $i, j, k$, we have
$$
s_{\mu, \lambda_{i,j}} =\mu \lambda_{i,j} \mu^{-1} \quad \textrm{and} \quad s_{\mu,\gamma_{k} } = 1. 
$$
Thus, the group $H_n$ is generated by the set
$$\{\, \mu \lambda_{i,j} \mu^{-1}  \mid \mu \in \mathbb{Z}_2^n~\textrm{and}~\ 1 \leq i <j \leq n \,\}. $$
To determine the defining relations of $H_n$, we apply the rewriting process to all the defining relations of $PTVT_n$ as given by Theorem \ref{sTVP_n}.

\begin{itemize}
\item Clearly, the relations $\gamma_k^2 =1$ and $\gamma_i\gamma_j \gamma_i \gamma_j=1$ in $PTVT_n$ give rise to the trivial relation in $H_n$.

\item Applying the rewriting process to the relation $\lambda_{i,j}\gamma_k\lambda_{i,j}^{-1}\gamma_k=1$ in $PTVT_n$, where $i,j,k$ are distinct indices, we get
\begin{eqnarray*}
\tau(\mu \lambda_{i,j}\gamma_k\lambda_{i,j}^{-1}\gamma_k\mu^{-1})
&=& s_{\overline{\mu}, \lambda_{i,j}} s_{\overline{\mu\lambda_{i,j}}, \gamma_{k}}
   s_{\overline{\mu\lambda_{i,j}\gamma_{k} \lambda_{i,j}^{-1}}, \lambda_{i,j}}^{-1}
   s_{\overline{\mu\lambda_{i,j}\gamma_{k}\lambda_{i,j}^{-1}}, \gamma_{k}} \\
&=& s_{\mu, \lambda_{i,j}} s_{\mu\gamma_{k}, \lambda_{i,j}}^{-1} \\
&=& (\mu \lambda_{i,j} \mu^{-1}) (\mu\gamma_{k} \lambda_{i,j} \gamma_{k} \mu^{-1})^{-1}.
\end{eqnarray*}
This yields the relation
$$
\lambda_{i,j}  = \gamma_{k}\lambda_{i,j}\gamma_{k} =\lambda_{i,j}^{\gamma_{k}}
$$
in $H_n$, which simply identifies two of its generators. Thus, $H_n$ is generated by the set
$$
\{\,  \lambda_{i,j},~~\lambda_{i,j}^{\gamma_i},~~  \lambda_{i,j}^{\gamma_j},~~ \lambda_{i,j}^{\gamma_i\gamma_j} \mid 1\leq i<j\leq n\,\}.
$$

\item  Applying the rewriting process to the relation $\lambda_{i,j}\gamma_i\gamma_j\lambda_{i,j}\gamma_j\gamma_i = 1$ in $PTVT_n$ gives
\begin{eqnarray*}
\tau(\mu \lambda_{i,j}\gamma_i\gamma_j\lambda_{i,j}\gamma_j\gamma_i\mu^{-1})
&=& s_{\overline{\mu}, \lambda_{i,j}} s_{\overline{\mu\lambda_{i,j}}, \gamma_{i}}
   s_{\overline{\mu\lambda_{i,j}\gamma_{i}}, \gamma_{j}} s_{\overline{\mu\lambda_{i,j}\gamma_{i}\gamma_j}, \lambda_{i,j}}
   s_{\overline{\mu\lambda_{i,j}\gamma_{i}\gamma_j\lambda_{i,j}}, \gamma_{j}}
   s_{\overline{\mu\lambda_{i,j}\gamma_{i}\gamma_j\lambda_{i,j}\gamma_j}, \gamma_{i}}\\
&=& s_{\mu, \lambda_{i,j}} s_{\mu\gamma_{i}\gamma_j, \lambda_{i,j}} \\
&=& (\mu \lambda_{i,j} \mu^{-1}) (\mu\gamma_{i}\gamma_j \lambda_{i,j} \gamma_j \gamma_{i}\mu^{-1}).
\end{eqnarray*}
This gives the relation
$$
 \lambda_{i,j}^{-1}  = \gamma_i\gamma_j\lambda_{i,j}\gamma_j\gamma_i=\lambda_{i,j}^{\gamma_i\gamma_j}
$$
in $H_n$, which identifies one generator with the inverse of the other. This implies that $H_n$ is generated by the smaller set
\begin{equation}\label{small gen set}
\{\,  \lambda_{i,j},~~\lambda_{i,j}^{\gamma_i},~~  \lambda_{i,j}^{\gamma_j} \mid 1\leq i<j\leq n\,\}.
\end{equation}

\item Finally, applying the rewriting process to the relation $\lambda_{i,j}\lambda_{k,l}\lambda_{i,j}^{-1}\lambda_{k,l}^{-1} = 1$ in $PTVT_n$, where $i,j,k, l$ are distinct indices,  we obtain
\begin{eqnarray*}
\tau(\mu \lambda_{i,j}\lambda_{k,l}\lambda_{i,j}^{-1}\lambda_{k,l}^{-1}\mu^{-1})
&=& s_{\overline{\mu}, \lambda_{i,j} }s_{\overline{\mu\lambda_{i,j}}, \lambda_{k,l}}
   s_{\overline{\mu\lambda_{i,j}\lambda_{k,l}\lambda_{i,j}^{-1}}, \lambda_{i,j}}^{-1}
   s_{\overline{\mu\lambda_{i,j}\lambda_{k,l}\lambda_{i,j}^{-1} \lambda_{k,l}^{-1}}, \lambda_{k,l}}^{-1} \\
&=& s_{\mu, \lambda_{i,j}}s_{\mu, \lambda_{k,l}}s_{\mu, \lambda_{i,j}}^{-1} s_{\mu, \lambda_{k,l}}^{-1} \\
&=& (\mu\lambda_{i,j} \mu^{-1}) (\mu \lambda_{k,l} \mu^{-1}) (\mu \lambda_{i,j} \mu^{-1})^{-1} ( \mu \lambda_{k,l} \mu^{-1})^{-1}.
\end{eqnarray*}
This gives the relation
$$
\lambda_{i,j}^{\mu} \lambda_{k,l}^{\mu} =\lambda_{k,l}^{\mu}\lambda_{i,j}^{\mu}
$$
in $H_n$, which are the desired defining relations among the generators from the set \eqref{small gen set}.
\end{itemize}
The presentation of  $PTVT_n$ obtained above shows that it is a right-angled Artin group. It is known that a right-angled Artin group is irreducible if and only if its defining graph does not decompose as a non-trivial simplicial join. The irreducibility of $PTVT_n$ follows from the observation that the complement of its defining graph is connected.
\end{proof}

\color{black}

\begin{corollary}\label{cor tvtn decomposition}
For each $n \ge 2$, the following decompositions hold:
    \begin{enumerate}
        \item $TVT_n \cong PTVT_n \rtimes S_n \cong  (H_n \rtimes \mathbb{Z}_2^n) \rtimes S_n$.
        \item $TVT_n \cong PTVT_n \rtimes S_n \cong  \Big(F_\infty * \mathbb{Z}_2\Big) \rtimes
\Big(F_\infty * \mathbb{Z}_2\Big) \rtimes \cdots \rtimes \Big(F_\infty * \mathbb{Z}_2\Big) \rtimes 
\big(\mathbb{Z} * \mathbb{Z}_2\big) \rtimes \mathbb{Z}_2 \rtimes S_n$, where $F_\infty$ is a free group of countably infinite rank.
    \end{enumerate}
\end{corollary}

\begin{corollary}\label{cor trivial center}
For each $n \ge 2$, we have $Z(TVT_n)=Z(PTVT_n)=1$.
\end{corollary}

\begin{proof}
Note that $PTVT_2 \cong \langle a,b,c \mid b^2=c^2=1,\; bc=cb,\; a^{-1}=bcabc\rangle$. Consider the subgroups 
$N=\langle a,b\rangle\cong \mathbb{Z}*\mathbb{Z}_2$ and $H=\langle c\rangle\cong \mathbb{Z}_2$ of $PTVT_2$. Then $PTVT_2 \cong N \rtimes H$ and
$Z(N)$ is trivial. If $z\in Z(PTVT_2)$, then we can write $z=n c^\varepsilon$ for some $n\in N$ and $\varepsilon\in\{0,1\}$.
\par

If $z=n$, then  $z\in Z(N)$, and hence $z=1$. If $z=n c$, then $n c x =z x =x z = x n c$  for all $x\in N$. Taking $x=b$ and using the fact that $b c = cb$, we get $n b = b n$. This implies that $ n \in \{1,b\}$.
\par

If $n=1$, then $cx=xc$. Taking $x=a$, the relation $a^{-1}=bcabc$ becomes $b a^{-1} b = a$. But, this is not possible in the free product $\mathbb{Z}*\mathbb{Z}_2$.
\par 
If $n=b$, then taking $x=a$, we obtain $bca = a bc$. In this case, the relation $a^{-1}=bcabc$ becomes $a^{-1}=b^2 a b^2$, which is not possible in $\mathbb{Z}*\mathbb{Z}_2$.  Hence, we must have $Z(PTVT_2)=1$.
\par

Recall that $Z(G) \subseteq Z(N)$ whenever $G = N \rtimes H$ and $Z(H) = 1$. Using Theorem~\ref{decom1} and the fact that $Z(F_\infty^{(n)} * \mathbb{Z}_2) =1$ for each $n \ge 3$, we obtain $Z(PTVT_n) =1$ for each $n \geq 2$. Similarly, since $Z(S_n) = 1$ for each $n \geq 2$,  in view of Corollary \ref{cor tvtn decomposition}, we obtain $Z(TVT_n) = 1$ for each $n \ge 2$.
\end{proof}

\begin{corollary}\label{cor:3.5}
For each $n \geq 2$, the groups $TVT_n$ and $PTVT_n$ are residually finite and Hopfian.
\end{corollary}

\begin{proof}
It is known that a right-angled Artin group is linear~\cite[Corollary~3.6]{MR1704150} and that a finitely generated linear group is residually finite~\cite{MR0003420}. Hence, $H_n$ is residually finite. Since any extension of a residually finite group by a finite group is again residually finite, it follows that $PTVT_n$ is residually finite. Further, $PTVT_n$ is Hopfian follows from the fact that every finitely generated residually finite group is Hopfian~\cite{MR0003420}. Similar reasoning implies that $TVT_n$ is residually finite and Hopfian for each $n \geq 2$.
\end{proof}
\medskip

\section*{Acknowledgements}
Singh is supported by the Swarna Jayanti Fellowship grants DST/SJF/MSA-02/2018-19 and SB/SJF/2019-20/04. Negi acknowledges Research Associateship from the grant SB/SJF/2019-20/04.

\end{document}